\documentclass{article}

\usepackage{amsmath,amssymb,amsfonts,esint}

\newcommand\dps{\displaystyle }

 \usepackage{latexsym}
 \usepackage{pdflscape}
 \usepackage{examplep}
 \usepackage{longtable}
 \usepackage{graphicx}
 \usepackage{amsmath}
 \usepackage{booktabs}
 \usepackage{memhfixc}
 \usepackage{dsfont}
 \usepackage{amsmath,amssymb,amsfonts,amsthm}
 \usepackage[francais,english]{babel} 
 \newtheorem{theorem}{Theorem}[section]
 \newtheorem{proposition}{Proposition}[section]
 
 \newtheorem{lemma}{Lemma}[section]
 \newtheorem{corollary}{Corollary}[section]
 \newtheorem{definition}{Definition}[section]
 \newtheorem{example}{Example}[section]

\def\cC{{\mathcal C}}
\def\N{{\mathbb N}}
\def\Z{{\mathbb Z}}
\def\R{{\mathbb R}}
\def\C{{\mathbb C}}
\def\QA{{Q(A)}}
\def\DA{{D(A)}}
\def\cH{{\cal H}}
\def\cR{{\cal R}}
\def\cC{{\cal C}}
\def\cT{{\cal T}}
\def\cL{{\cal L}}
\def\cS{{\cal S}}
\def\cP{{\cal P}}
\def\cA{{\cal A}}
\def\cM{{\cal M}}
\def\cZ{{\cal Z}}
\def\cG{{\cal G}}
\def\cI{{\cal I}}
\def\ah{\widehat{a}}
\def\uh{\widehat{u}}
\def\mh{\widehat{m}}
\def\at{\widetilde{a}}
\def\mt{\widetilde{m}}
\def\wlim{\rightharpoonup}
\def\sigess{\widehat{\sigma}_{\rm ess}}

\def\cD{\mathcal{D}}
\def\fP{\mathfrak{P}}
\def\un{\mathds{1}}

\title{Non-consistent approximations of self-adjoint eigenproblems: \\
Application to the supercell method\thanks{This work was financially supported by the ANR grant MANIF.}}

\author{Eric Canc\`es\thanks{Universit\'e Paris Est, CERMICS, Projet MICMAC, Ecole des Ponts ParisTech - INRIA, 6 \& 8 avenue Blaise Pascal, 77455 Marne-la-Vall\`ee Cedex 2, France, ({\tt cances@cermics.enpc.fr}, {\tt ehrlachv@cermics.enpc.fr})} \and Virginie Ehrlacher \and Yvon Maday\thanks{Universit\'e Pierre et Marie Curie-Paris 6, UMR 7598, Laboratoire J.-L. Lions, Paris, F-75005 France, and Division of Applied Mathematics, Brown University,
182 George Street, Providence, RI 02912, USA, ({\tt maday@ann.jussieu.fr})}}

\begin{document}

\maketitle

\begin{abstract} 
In this article, we introduce a general theoretical framework to analyze non-consistent approximations of the discrete eigenmodes of a self-adjoint operator. We focus in particular on the discrete eigenvalues laying in spectral gaps. We first provide {\it a priori} error estimates on the eigenvalues and eigenvectors in the absence of spectral pollution. We then show that the supercell method for perturbed periodic Schr\"odinger operators falls into the scope of our study. We prove that this method is spectral pollution free, and we derive optimal convergence rates for the planewave discretization method, taking numerical integration errors into account. Some numerical illustrations are provided.
\end{abstract}

%\begin{keywords} 
%...
%\end{keywords}

%\begin{AMS}
%...
%\end{AMS}

%\pagestyle{myheadings}
%\thispagestyle{plain}
%\markboth{E. Canc\`es, V. Ehrlacher and Y. Maday}{Numerical analysis of perturbed periodic Schr\"odinger operators}

\section{Introduction}

This article is concerned with the numerical analysis of the computation of the discrete eigenmodes of a self-adjoint operator $A$, on an infinite dimensional separable Hilbert space $\cH$. The focus is particularly set on the eigenmodes corresponding to discrete eigenvalues located in spectral gaps. 

The main application we have in mind is concerned with perturbed periodic Schr\"odinger operators of the form 
$$ 
A:= -\Delta + V_{\rm per} +W,
$$
where $\Delta$ is the Laplace operator on $L^2(\R^d)$, $V_{\rm per}$ a periodic function of $L^p_{\rm loc}(\R^d)$ with $p=2$ if $d\leq 3$, $p>2$ for $d=4$ and $p=d/2$ for $d\geq 5$, and $W\in L^{\infty}(\R^d)$ a perturbation of the potential going to zero at infinity. The operator $A$ is self-adjoint and bounded from below on $\cH:= L^2(\R^d)$ with domain $H^2(\R^d)$. Perturbed periodic Schr\"odinger operators are encountered in electronic structure theory, and in the study of photonic crystals. In the case of a perfectly periodic crystal ($W=0$), the spectrum of the operator $A^0:= - \Delta + V_{\rm per}$ is purely absolutely continuous, and composed of a union 
of intervals of $\R$. It follows from Weyl's theorem~\cite{ReedSimon4} that the essential spectra of $A$ and $A^0$ are identical. On the other hand, when $W \neq 0$, some discrete eigenvalues may appear in the band gaps of the spectrum of $A$. The corresponding eigenmodes, which can be interpreted as bound states trapped by local defects, are difficult to compute for numerical methods can produce spectral pollution.   

In a general theoretical framework, the eigenvalues of $A$ and the associated eigenvectors can be obtained by solving the variational problem
\begin{equation}\label{eq:cont}
\left\{
\begin{array}{l}
 \mbox{\rm find }(\psi, \lambda) \in \QA\times \R \mbox{ such that}\\
\forall \phi\in \QA, \; a(\psi, \phi) = \lambda m(\psi, \phi),\\
\end{array}
\right .
\end{equation}
where $m(\cdot, \cdot)$ is the scalar product of $\cH$, $\QA$ the form domain of $A$, and $a(\cdot, \cdot)$ the sesquilinear form associated with $A$ (see for instance \cite{Chatelin}).

A sequence $(X_n)_{n\in\N}$ of finite dimensional approximation subspaces of $\QA$ being given, we consider for all $n\in\N$, the self-adjoint operator $A|_{X_n}: X_n \to X_n$ defined by
$$
\forall (\psi_n, \phi_n) \in X_n \times X_n, \quad m\left( A|_{X_n} \psi_n, \phi_n\right) = a(\psi_n, \phi_n).
$$
The standard Galerkin method consists in approximating the discrete eigenvalues of the operator $A$ by the eigenvalues of the discretized operators $A|_{X_n}$, 
the latter being obtained by solving the variational problem
$$
\left\{ 
\begin{array}{l}
 \mbox{\rm find } (\psi_n, \lambda_n)\in X_n \times \R \mbox{ \rm such that}\\
\forall \phi_n \in X_n, \; a(\psi_n, \phi_n) = \lambda_n m(\psi_n, \phi_n).\\ 
\end{array}
\right .
$$
 According to the Rayleigh-Ritz theorem \cite{ReedSimon4}, under the natural assumption that the sequence $(X_n)_{n\in\N}$ satisfies 
$$
\forall \phi\in \QA, \quad \mathop{\inf}_{\phi_n\in X_n} \|\phi - \phi_n\|_{\QA} \mathop{\longrightarrow}_{n\to\infty} 0, 
$$
this method allows to compute the eigenmodes of $A$ associated with the discrete eigenvalues located below the bottom of the essential spectrum. 
Indeed, if $A$ is bounded below and possesses exactly $M$ discrete eigenvalues $\lambda_1 \leq \lambda_2 \leq \cdots \leq \lambda_M$ 
(taking multiplicities into account) lower than $\min \sigma_{\rm ess}(A)$, where $\sigma_{\rm ess}(A)$ denotes the essential spectrum of $A$, and if $\left\{ \lambda_j^n \right\}_{1\leq j \leq \mbox{\rm dim }X_n }$ are the eigenvalues of $A|_{X_n}$, it is well-known that
$$
\forall 1\leq j \leq M, \; \lambda_j^n \mathop{\downarrow}_{n\to\infty} \lambda_j.
$$
The situation is much more delicate when one tries to approximate eigenvalues which are located in spectral gaps of $A$ since 
$$
\forall M <j \leq \mbox{\rm dim }X_n, \quad \lambda_j^n \mathop{\downarrow}_{n\to\infty} \min \sigma_{\rm ess}(A).
$$
When dealing with the approximation of discrete eigenvalues of $A$ located in spectral gaps, 
the standard Galerkin method may give rise to spectral pollution: some sequences $(\lambda_n)_{n \in \N}$, where for each
 $n$, $\lambda_n \in \sigma(A|_{X_n})$, may converge to real numbers which do not belong to the spectrum of $A$. Spectral pollution occurs in a broad variety of physical settings, including elasticity theory, electromagnetism, hydrodynamics and quantum physics \cite{Arnold, Boffi, Dauge, Rapaz,Shabaev}, and has been extensively studied in the framework of the standard Galerkin method~\cite{Boulton07,BoultonBoussaid,BBL, Davies, Descloux,Rappaz1, Rappaz2, Hansen,LevSha}. We refer to \cite{BoultonLevitin,CEM,LS} for an analysis of spectral pollution for perturbed periodic Schr\"odinger operators. 

On the other hand, few results have been published on the numerical computation of eigenmodes in spectral gaps by means 
of non-consistent methods, based on generalized eigenvalue problems of the form
$$
 \left\{
\begin{array}{l}
 \mbox{\rm find }(\psi_n, \lambda_n)\in X_n\times \R \mbox{ \rm such that}\\
\forall \phi_n\in X_n, \quad a_n(\psi_n, \phi_n) = \lambda_n m_n(\psi_n, \phi_n),\\
\end{array}
\right .
$$
where for all $n\in\N$, $a_n(\cdot, \cdot)$ and $m_n(\cdot, \cdot)$ are symmetric bilinear forms on $\QA$, {\it a priori} different from $a(\cdot, \cdot)$ and $m(\cdot, \cdot)$. 

\medskip

In this article, we consider a general theoretical framework to analyze non-consistent methods for the computation of the discrete eigenmodes of a self-adjoint operator. 
After introducing some notation and definitions in Sections~\ref{sec:not} and~\ref{sec:approximation}, we state our main result (Theorem~\ref{th:thmain}) in Section~\ref{sec:theory}. Theorem~\ref{th:thmain} provides {\it a priori} error estimates on the eigenvalues and eigenvectors in the absence of spectral pollution. Its proof is given in Section~\ref{sec:proofmain}. 

In Section~\ref{sec:supercell}, we show that the supercell method for perturbed periodic Schr\"odinger operators falls into the scope of Theorem~\ref{th:thmain}. We prove that this method is spectral pollution free, and we derive optimal convergence rates for the planewave discretization method, taking numerical integration errors into account. The corresponding proofs are detailed in Section~\ref{sec:supercellproof}, and some numerical illustrations are provided in Section~\ref{sec:numerical}.

\section{Approximations of a self-adjoint operator}

\subsection{Some notation}\label{sec:not}

Throughout this paper, $\cH$ denotes a separable Hilbert space, endowed with the scalar product $m(\cdot, \cdot)$ and associated norm $\|\cdot\|_{\cH}$, and $A$
 a self-adjoint operator on $\cH$ with dense domain $\DA$. We denote by $\QA:=D(|A|^{1/2})$ the form domain of $A$ and by $a(\cdot,\cdot)$ the symmetric bilinear form on 
$\QA$ associated with $A$. Recall that the vector space $\QA$, endowed with the scalar product $\langle \cdot, \cdot \rangle_{\QA}$, defined as
$$
\forall \psi,\phi\in\QA, \quad \langle \psi,\phi\rangle_{\QA} := m(\psi,\phi) + m \left(|A|^{1/2} \psi, |A|^{1/2} \phi\right),
$$
is a Hilbert space; the associated norm is denoted by $\|\cdot\|_{\QA}$.

\medskip

\begin{example}
Perturbed periodic Schr\"odinger operators $A:= -\Delta + V_{\rm per} + W$ are self-adjoint semibounded operators on $\cH:=L^2(\R^d)$, with domain $\DA:= H^2(\R^d)$ and form domain 
$\QA := H^1(\R^d)$.
\end{example}

\medskip

For any finite dimensional vector subspace $X$ of $\cH$ such that $X \subset \QA$, we introduce the following notation
\begin{itemize}
\item $i_X \, : \, X \hookrightarrow \cH$ is the canonical embedding of $X$ into $\cH$;
\item $i^*_X \, : \, \cH \rightarrow X$ is the adjoint of $i_X$, that is the orthogonal projection from $\cH$ onto $X$ associated with the scalar product $m\left( \cdot, \cdot \right)$;
\item $A|_{X} \, : \, X \rightarrow X$ is the self-adjoint operator on $X$ defined by 
$$
\forall (\psi,\phi) \in X \times X, \quad m\left( A|_{X}\psi,\phi\right)=a(\psi,\phi);
$$ 
\item $\Pi_X^{\cH} \, : \, \cH \rightarrow \cH$ and $\Pi_X^{\QA}: \, \QA \rightarrow \QA$ are the orthogonal projections onto $X$ for
 $\left(\cH, m( \cdot, \cdot )\right)$ and $\left(\QA, \langle \cdot, \cdot\rangle_{\QA}\right)$, respectively. 
\end{itemize}

We set
$$
\sigess(A) := \overline{\sigma(A)}^{\overline{\R}} \setminus \sigma_{\rm d}(A), 
$$
where $\sigma_{\rm d}(A)$ is the discrete spectrum 
of $A$, and where $\overline{\sigma(A)}^{\overline{\R}}$ is the closure of $\sigma(A)$, the spectrum of $A$, in $\overline\R:=\R \cup \left\{\pm\infty\right\}$. 
A {\it spectral gap} of $A$ is an interval $\left(\Sigma^-, \Sigma^+\right)$ such that
 $\Sigma^-,\Sigma^+ \in \sigess(A) \cap \R$ and $\left(\Sigma^-, \Sigma^+\right) \cap \sigess(A) = \emptyset$ (which implies that $\mbox{\rm Tr}(\un_{(-\infty, \Sigma^-]}(A)) = \mbox{\rm Tr}(\un_{[\Sigma^+, \infty)}(A)) = \infty$). As usual, $\un_B$ denotes the characteristic function of the Borel set $B\subset \R$.
The discrete eigenvalues of the operator $A$ in a spectral gap $\left(\Sigma^-, \Sigma^+\right)$, if any, are isolated and of finite multiplicities,
 but can accumulate at $\Sigma^-$ and/or $\Sigma^+$ \cite{ReedSimon4}. 

\medskip

Let us finally recall the notions of limit superior and limit inferior of a sequence of sets of complex numbers (see for instance \cite{Chatelin}).

\begin{definition}
Let $(E_n)_{n\in\N}$ be a sequence of subsets of $\C$.
\begin{itemize}
 \item The set $\dps \mathop{\overline{\lim}}_{n\to\infty} E_n$ (limit superior) is the set of all complex numbers $\lambda\in \C$ such that there exist a subsequence $(E_{n_k})_{k\in\N}$ of 
$(E_n)_{n\in\N}$ and a sequence $(\lambda_{n_k})_{k\in\N}$ of complex numbers such that for all $k\in\N$, $\lambda_{n_k} \in E_{n_k}$ and 
$\dps \mathop{\lim}_{k\to\infty} \lambda_{n_k} = \lambda$. 
\item The set  $\dps \mathop{\underline{\lim}}_{n\to\infty} E_n$ (limit inferior) is the set of all complex numbers $\lambda\in \C$ such that there exists 
a sequence $(\lambda_{n})_{n\in\N}$ of complex numbers such that for all $n\in\N$, $\lambda_{n} \in E_{n}$ and 
$\dps \mathop{\lim}_{n\to\infty} \lambda_{n} = \lambda$.
\item If $\dps \mathop{\underline{\lim}}_{n\to\infty} E_n = \mathop{\overline{\lim}}_{n\to\infty} E_n$, then 
$
\dps \mathop{\lim}_{n\to\infty} E_n:=\mathop{\underline{\lim}}_{n\to\infty} E_n = \mathop{\overline{\lim}}_{n\to\infty} E_n.
$
\end{itemize}
\end{definition}

\subsection{Consistent and non-consistent approximations}\label{sec:approximation}

\begin{definition}
An approximation $(\cT_n)_{n\in\N}$ of a self-adjoint operator $A$ is a sequence such that, for all $n\in\N$,
$$
\cT_n: = (X_n,a_n,m_n),
$$
where
\begin{itemize}
 \item $(X_n)_{n\in\N}$ is a sequence of finite dimensional subspaces of $\QA$;
\item $(a_n)_{n\in\N}$ is a sequence of symmetric bilinear forms on $\QA$;
\item $(m_n)_{n\in\N}$ is a sequence of symmetric bilinear forms on $\QA$ such that the restriction of $m_n$ to $X_n$ forms a scalar product on $X_n$. We denote by $\|\cdot\|_{X_n}$ the associated norm: $\forall \phi_n \in X_n$, $\|\phi_n\|_{X_n} = m_n(\phi_n,\phi_n)^{1/2}$.
\end{itemize}
The approximation $(\cT_n)_{n\in\N}$ is called consistent if, for any $(\psi,\lambda)$ solution of (\ref{eq:cont}),
$$
\forall \phi_n \in X_n, \quad a_n(\psi, \phi_n) = \lambda m_n(\psi, \phi_n),
$$
and non-consistent otherwise.
\end{definition}

The approximation $(\cT_n)_{n\in\N}$ is referred to as a standard Galerkin method if, for all $n\in\N$, $a_n = a$ and $m_n = m$. Standard Galerkin methods are obviously consistent.

If $(\cT_n)_{n\in\N}$ is an approximation of $A$, we denote by $\cA_n$ and $\cM_n$ the $m$-symmetric (i.e. symmetric w.r.t. the scalar product $m(\cdot,\cdot)$) linear operators on $X_n$ defined by: $\forall \phi_n, \psi_n\in X_n$, 
\begin{eqnarray*}
&& m\left( \cA_n \phi_n, \psi_n \right) =  a_n(\phi_n, \psi_n),\\
&& m\left( \cM_n \phi_n, \psi_n \right) =  m_n(\phi_n, \psi_n).
\end{eqnarray*}
Since $m_n$ is a scalar product on $X_n$, the operator $\cM_n$ is invertible and we can define the operator
$$
A_n = \cM_n^{-1/2} \cA_n \cM_n^{-1/2}
$$
on $X_n$, which is $m$-symmetric as well. The generalized eigenvalue problem
\begin{equation}\label{eq:geneig}
\left\{
\begin{array}{l}
 \mbox{find }(\psi_n, \lambda_n)\in X_n \times \R \mbox{ such that }\|\psi_n\|_{X_n}^2 = 1 \mbox{ and }\\
\forall \phi_n \in X_n, \; a_n(\psi_n, \phi_n) = \lambda_n m_n(\psi_n, \phi_n),\\
\end{array}
 \right .
\end{equation}
is then equivalent, through the change of variable $\xi_n=\cM_n^{1/2}\psi_n$, to the eigenvalue problem
$$
\left\{
\begin{array}{l}
 \mbox{find }(\xi_n, \lambda_n)\in X_n \times \R \mbox{ such that }\|\xi_n\|_{\cH}^2 = 1 \mbox{ and }\\
A_n\xi_n = \lambda_n \xi_n.\\
\end{array}
\right.
$$

The main objective of this work is to provide sufficient conditions on such potentially non-consistent approximations $(\cT_n)_{n\in\N}$ so that the discrete eigenvalues of $A$ and the associated eigenvectors are well-approximated in a certain sense by eigenvalues and eigenvectors of the discretized problems (\ref{eq:geneig}). We wish to provide a framework which will enable us to deal with 
the supercell method for perturbed periodic linear Schr\"odinger operators described in Section~\ref{sec:supercell}.

\section{An abstract convergence result}\label{sec:theory}

\subsection{The general case}

Let us consider an approximation $(\cT_n)_{n\in\N}$ of $A$ satisfying the following assumptions:
\begin{itemize}
 \item [(A1)]  $\dps \forall \psi \in \QA, \; \left\| \left( 1 - \Pi_{X_n}^{\QA}\right) \psi \right\|_{\QA} \mathop{\longrightarrow}_{n\to\infty} 0$;
\item[(A2)] there exists $0 < \gamma \le \Gamma < \infty$ such that for all $n\in\N$ and all $\psi_n ,\phi_n \in X_n$, 
\begin{eqnarray*}
 \gamma \|\psi_n\|_{\cH}^2 \leq & m_n(\psi_n, \psi_n) & \leq \Gamma \|\psi_n\|_{\cH}^2,\\
& |a_n(\psi_n, \phi_n)| & \leq  \Gamma\|\psi_n\|_{\QA}\|\phi_n\|_{\QA};
\end{eqnarray*}
\item[(A3)] for any compact subset $K\subset\C$, if there exists a subsequence $(\cT_{n_k})_{k\in\N}$ of $(\cT_n)_{n\in\N}$ such that 
$\mbox{\rm dist }(K, \sigma(A_{n_k})) \geq \alpha_K$ for some $\alpha_K >0$ independent of $k\in\N$, then there exists $c_K >0$ such that for all $\mu\in K$ and all $k\in \N$, 
$$
\mathop{\inf}_{w_{n_k}\in X_{n_k}} \mathop{\sup}_{v_{n_k} \in X_{n_k}} \frac{|(a_{n_k} - \mu m_{n_k})(w_{n_k}, v_{n_k})|}{\|w_{n_k}\|_{\QA} \|v_{n_k}\|_{\QA}} \geq c_K;
$$  
\item[(A4)] there exist $\kappa \in \R_+$ and, for each $n\in\N$, two symmetric bilinear forms $\widetilde{a}_n$ and $\widetilde{m}_n$ on $\QA$, and four seminorms $r_n^a$, $r_n^m$, $s_n^a$, $s_n^m$ on $\QA$ such that $\forall \phi_n, \psi_n \in X_n$, 
$$
\begin{array}{rcl}
 \gamma \|\psi_n\|_{\cH}^2 \leq & \widetilde{m}_n(\psi_n, \psi_n) & \leq \Gamma \|\psi_n\|_{\cH}^2,\\
& |\widetilde{a}_n(\psi_n, \phi_n)| & \leq  \Gamma\|\psi_n\|_{\QA}\|\phi_n\|_{\QA},\\
\end{array}
$$
and $\forall \phi,\psi \in \QA$,
$$
\begin{array}{ll}
 \dps |(a-\widetilde{a}_n)(\phi, \psi)| \leq r_n^a(\phi)r_n^a(\psi), & \qquad \dps |(m- \widetilde{m}_n) (\phi, \psi)| \leq  r_n^m(\phi)r_n^m(\psi),\\
\dps r_n^a(\phi) \leq \kappa \|\phi\|_{\QA} , & \qquad \dps r_n^m(\phi) \leq \kappa \|\phi\|_{\cH},\\
\dps r_n^a(\phi) \mathop{\longrightarrow}_{n\to\infty} 0 , & \qquad \dps r_n^m(\phi) \mathop{\longrightarrow}_{n\to\infty} 0,\\
\end{array}
$$
and 
$$ \!\!\!\!\!\!\!\!\!\!\!\!\!\!\!\!\!
\begin{array}{ll}
\dps \mathop{\sup}_{w_n\in X_n} \frac{|(a_n-\widetilde{a}_n)(\Pi_{X_n}^\QA \phi, w_n)|}{\|w_n\|_{\QA}} \leq s_n^a(\phi), & \quad \dps \mathop{\sup}_{w_n\in X_n} \frac{|(m_n- \widetilde{m}_n)(\Pi_{X_n}^\QA\phi, w_n)|}{\|w_n\|_{\QA}} \leq  s_n^m(\phi),\\ \\
\dps s_n^a(\phi) \leq \kappa \|\phi\|_{\QA} , &\quad \dps s_n^m(\phi) \leq \kappa \|\phi\|_{\cH},\\
\dps s_n^a(\phi) \mathop{\longrightarrow}_{n\to\infty} 0 , &\quad \dps s_n^m(\phi) \mathop{\longrightarrow}_{n\to\infty} 0.\\
\end{array}
$$
\end{itemize}

Before stating our main result, let us comment on these assumptions.

\medskip

Conditions (A1) and (A2) are classical. The former means that any $\psi \in Q(A)$ can be approximated in $Q(A)$ by a sequence $(\psi_n)_{n \in \N}$ such that $\psi_n \in X_n$ for each $n \in \N$. The latter ensures that, {\em uniformly in $n$}, the norms $\|\cdot\|_{X_n}$ and $\|\cdot\|_{\cal H}$ are equivalent on $X_n$, and the bilinear forms $a_n$ are continuous on $X_n$, the space $X_n$ being endowed with the norm $\|\cdot\|_{Q(A)}$.

\medskip

Assumption (A3) is important in our proof since it enables us to apply Strang's lemma (see Section~\ref{sec:appendix}) with a uniform discrete inf-sup condition. 

For the supercell approximation, we will prove a stronger result:
\begin{itemize}
 \item [(A3')] for any compact subset $K \subset \C$, there exists $c_K>0$ such that for all $n\in\N$ and all $\mu\in K$, 
$$
\mathop{\inf}_{w_n\in X_n} \mathop{\sup}_{v_n\in X_n} \frac{|(a_n- \mu m_n)(w_n, v_n)|}{\|w_n\|_{\QA} \|v_n\|_{\QA}}\geq c_K \min(1, \mbox{\rm dist}(\mu, \sigma(A_n))).
$$
\end{itemize}
 It is easily checked that (A3') implies (A3).

\medskip

Let us finally comment on condition (A4) in the perspective of the analysis of the supercell method with numerical integration addressed in Section~\ref{sec:supercell}.
 In the latter setting, the introduction of the bilinear forms $\widetilde{a}_n$ and $\widetilde{m}_n$ aims at separating in the error bounds of Theorem~\ref{th:thmain} 
the contributions inherently due to the supercell method (truncation of the domain and artificial periodic boundary conditions) and those due to numerical integration.
 We postpone until Section~\ref{sec:supercell} the precise definitions of $a_n$, $m_n$, $\widetilde{a}_n$ and $\widetilde{m}_n$ in this context. 

Note that (A4) implies that the approximation $(\cT_n)_{n\in\N}$ is weakly consistent in the sense that for all $\phi\in\QA$, the consistency errors $r_n^a(\phi)$, $r_n^m(\phi)$, 
$s_n^a(\phi)$ and $s_n^m(\phi)$ converge to $0$ as $n$ goes to infinity.

\medskip

We are now in position to state our main result.

\begin{theorem}\label{th:thmain}
 Let $A$ be a self-adjoint operator on $\cH$, $\lambda \in \sigma_{\rm d}(A)$ a discrete eigenvalue of $A$ with multiplicity $q$, and $(\cT_n)_{n\in\N}$ an approximation of $A$ satisfying assumptions (A1)-(A4). 
Then,
\begin{itemize}
 \item[1.] \bfseries Convergence of the eigenvalues \normalfont \itshape
\end{itemize}
\begin{equation}\label{eq:atteinte}
\lambda \in \mathop{\underline{\lim}}\sigma(A_n).
\end{equation}
\begin{itemize}
\item[2.] \bfseries A priori error estimates in the absence of spectral pollution \normalfont \itshape 
\end{itemize}
Assume that
$$
\mbox{\rm (B1)} \quad \exists \varepsilon >0 \mbox{ {\rm s.t.} } \dps (\lambda-\varepsilon, \lambda +  \varepsilon) \cap \sigma(A) = \{ \lambda\} \mbox{ {\rm and} } \mathop{\overline{\lim}}_{n\to\infty} \sigma(A_n) \cap (\lambda - \varepsilon, \lambda + \varepsilon) = \{\lambda\}.
$$
Let $\cP: = \un_{\left\{\lambda\right\}}(A)$ be the orthogonal projection on $\mbox{\rm Ker}(A -\lambda)$ and 
$$\cP_n:= i_{X_n} \cM_n^{-1/2} \un_{(\lambda - \varepsilon/2, \lambda + \varepsilon/2)}(A_n)\cM_n^{1/2} i_{X_n}^*.
$$ 
Then,
\begin{equation}\label{eq:dim}
 \mbox{\rm Rank}(\cP_n) \geq q,
\end{equation}
and there exists $C\in \R_+$ such that, for $n$ large enough, 
\begin{equation}\label{eq:proj1}
 \|(\cP-\cP_n)\cP\|_{\cL(\cH, \QA)} \leq C \left(\left\| \left( 1 - \Pi_{X_n}^{\QA}\right) \cP \right\|_{\cL(\cH, \QA)} + \cR_n^a + \cR_n^m + \cS_n^a + \cS_n^m\right),
\end{equation}
with
\begin{eqnarray*}
 \cR_n^a & := & \dps \mathop{\sup}_{\psi \in {\rm Ran}(\cP), \; \|\psi\|_{\cH} = 1} r_n^a(\psi),\\
 \cR_n^m & := & \dps \mathop{\sup}_{\psi \in {\rm Ran}(\cP), \; \|\psi\|_{\cH} = 1} r_n^m(\psi),\\
 \cS_n^a & := & \dps \mathop{\sup}_{\psi \in {\rm Ran}(\cP), \; \|\psi\|_{\cH} = 1} s_n^a(\psi),\\
 \cS_n^m & := & \dps \mathop{\sup}_{\psi \in {\rm Ran}(\cP), \; \|\psi\|_{\cH} = 1} s_n^m(\psi).
\end{eqnarray*}
If we assume in addition that
$$  
\mbox{\rm (B2)} \quad \mbox{for $n$ large enough, $\mbox{\rm Rank}(\cP_n) = q$},
$$
then there exists $C \in \R_+$ such that, for $n$ large enough, 
\begin{equation}\label{eq:proj2}
 \|(\cP-\cP_n)\cP_n\|_{\cL(\cH,\QA)} \leq  C \left(\left\| \left( 1 - \Pi_{X_n}^{\QA}\right) \cP \right\|_{\cL(\cH, \QA)} + \cR_n^a + \cR_n^m + \cS_n^a + \cS_n^m\right),
\end{equation}
\begin{equation}\label{eq:eigenvalue}
 \dps \mathop{\max}_{\lambda_n \in \sigma(A_n)\cap (\lambda - \varepsilon/2, \lambda+\varepsilon/2)} |\lambda_n - \lambda| \leq C \left( \left(\left\| \left( 1 - \Pi_{X_n}^{\QA} \right) \cP \right\|_{\cL(\cH, \QA)} + \cR_n^a + \cR_n^m \right)^2 + \cS_n^a + \cS_n^m \right). 
\end{equation}
\end{theorem}

\medskip

It is easy to check that $P_n := \cM_n^{-1/2} \un_{(\lambda - \varepsilon/2, \lambda + \varepsilon/2)}(A_n)\cM_n^{1/2}$ is the $m_n$-orthogonal projection of $X_n$ onto the space $Y_n \subset X_n$ spanned by the eigenvectors of (\ref{eq:geneig}) associated with the eigenvalues belonging to the interval $(\lambda - \varepsilon/2, \lambda +\varepsilon/2)$. The operator $\cP_n = i_{X_n} P_n i_{X_n}^\ast \in {\cal L}(\cH)$ is therefore a (non-orthogonal) projection on the finite dimensional space $i_{X_n} Y_n \subset \cH$. 

\medskip

Theorem~\ref{th:thmain} implies that, if $(\cT_n)_{n\in\N}$ is an approximation of the operator $A$ satisfying (A1)-(A4), for all discrete eigenvalue $\lambda$ of $A$, there exists a sequence $(\lambda_n)_{n\in\N}$ of elements of $\sigma(A_n)$ converging to $\lambda$. Assumption (B1) states that there is no spurious eigenvalues in the vicinity of~$\lambda$.
Estimate (\ref{eq:proj1}) shows that under assumption (B1), for each eigenvector $\psi$ of $A$ associated with the discrete eigenvalue $\lambda$, there exists a sequence 
$(\psi_n)_{n\in\N}$ of elements of $\mbox{\rm Ran}(\cP_n)$ which strongly converges towards $\psi$ in $\QA$.

\medskip

On the other hand, there may {\it a priori} exist a sequence $(\psi_n)_{n\in\N}$ of normalized elements of $\mbox{\rm Ran}(\cP_n)$ weakly converging in $\cH$ towards a vector 
that is not an eigenvector of $A$ associated with~$\lambda$. This is excluded when we make the additional assumption (B2). Assumption (B2) means that, for $n$ large enough, the sum of the multiplicities of the eigenvalues of $A_n$ close to $\lambda$ is equal to the multiplicity $q$ of $\lambda$. Under this assumption, if $(\psi_n)_{n\in\N}$ is a sequence of 
vectors of $\cH$ such that for each $n$ large enough, $\psi_n$ is an $\cH$-normalized eigenvector of $A_n$ associated with an eigenvalue $\lambda_n\in (\lambda - \varepsilon, \lambda+\varepsilon)$, and if $(\psi_n)_{n \in \N}$ weakly converges in $\cH$ towards some $\psi\in\cH$, then estimate (\ref{eq:proj2}) implies that $\psi$ is a $\cH$-normalized eigenvector of $A$ associated with the eigenvalue $\lambda$ 
and that the convergence of $(\psi_n)_{n\in\N}$ to $\psi$ holds strongly in $\QA$. 

\medskip

Lastly, estimate (\ref{eq:eigenvalue}) shows that when $\widetilde{a}_n=a_n$ and $\widetilde{m}_n=m_n$ (which is the case in the supercell model when numerical integration errors are neglected), then $\cS_n^a = \cS_n^m = 0$, and the convergence rate of the eigenvalues is twice the convergence rate of the eigenvectors measured in the $\QA$ norm. Such a doubling of the convergence rate is expected in variational approximations of linear eigenvalue problems~(see e.g.~\cite{Chatelin}).

\subsection{Standard Galerkin method}

Let us now consider the special case when for all $n\in\N$, 
$\cT_n = (X_n, a, m)$ where $(X_n)_{n\in\N}$ is a sequence of finite dimensional subspaces of $\QA$ satisfying (A1). In this case, for all $n\in\N$, $A_n = A|_{X_n}$, ${\cal M}_n$ is the identity operator, and
$$
\cP_n= i_{X_n}\un_{(\lambda - \varepsilon/2, \lambda + \varepsilon/2)}(A|_{X_n}) i_{X_n}^* 
$$
is an orthogonal projector with respect to the scalar product $m$. 

\medskip

In this setting, (A2) and (A4) are obviously satisfied, and (A3) and (A3') respectively read
\begin{itemize}
 \item [(C3)]  for any compact subset $K\subset\C$, if there exists a subsequence $(X_{n_k})_{k\in\N}$ of $(X_n)_{n\in\N}$ such that 
$\mbox{\rm dist }(K, \sigma(A|_{X_{n_k}})) \geq \alpha_K$ for some $\alpha_K >0$ independent of $k\in\N$, then there exists $c_K >0$ such that for all $\mu\in K$ and all $k\in \N$, 
$$
\mathop{\inf}_{w_{n_k}\in X_{n_k}} \mathop{\sup}_{v_{n_k} \in X_{n_k}} \frac{|(a - \mu m)(w_{n_k}, v_{n_k})|}{\|w_{n_k}\|_{\QA} \|v_{n_k}\|_{\QA}} \geq c_K;
$$  
\end{itemize}
and
\begin{itemize}
 \item [(C3')] for all  compact subset $K \subset \C$, there exists $c_K>0$ such that for all $n\in\N$ and all $\mu\in K$, 
$$
\mathop{\inf}_{w_n\in X_n} \mathop{\sup}_{v_n\in X_n} \frac{|(a- \mu m)(w_n, v_n)|}{\|w_n\|_{\QA} \|v_n\|_{\QA}}\geq c_K \min(1, \mbox{\rm dist}(\mu, \sigma(A|_{X_n}))).
$$
\end{itemize}
It is proved in \cite{these} that, when $A$ is semibounded, (C3'), and thus (C3), automatically hold. 
On the other hand, when $A$ is not semibounded, (C3) is not always satisfied. An explicit counterexample is given in \cite{these}.

\medskip

The formulation of Theorem~\ref{th:thmain} simplifies in this case as follows:
\begin{corollary}\label{th:thexact}
 Let $A$ be a self-adjoint operator on $\cH$, $\lambda \in \sigma_{\rm d}(A)$ a discrete eigenvalue of $A$ with multiplicity $q$, and $(X_n)_{n\in\N}$ a sequence of finite dimensional subspaces of $\QA$ such that
$$ 
\forall \psi \in \QA, \quad \left\| \left( 1 - \Pi_{X_n}^{\QA}\right) \psi \right\|_{\QA} \mathop{\longrightarrow}_{n\to\infty} 0.
$$ 
Let us assume that either $A$ is semibounded or $(X_n)_{n\in\N}$ satisfies assumption (C3). 
Then,
\begin{itemize}
 \item[1.] \bfseries Convergence of the eigenvalues \normalfont \itshape
\end{itemize}
$$
\lambda \in \mathop{\underline{\lim}}\sigma(A|_{X_n}).
$$
\begin{itemize}
\item[2.] \bfseries A priori error estimates in the absence of spectral pollution \normalfont \itshape 
\end{itemize}
Assume that
$$
\mbox{\rm (D1)} \quad \exists \varepsilon >0 \mbox{ {\rm s.t.} } \dps (\lambda-\varepsilon, \lambda +  \varepsilon) \cap \sigma(A) = \{ \lambda\} \mbox{ {\rm and} } \mathop{\overline{\lim}}_{n\to\infty} \sigma(A_n) \cap (\lambda - \varepsilon, \lambda + \varepsilon) = \{\lambda\}.
$$
Let $\cP: = \un_{\left\{\lambda\right\}}(A)$ be the orthogonal projection on $\mbox{\rm Ker}(A -\lambda)$ and $\cP_n:= i_{X_n} \un_{(\lambda - \varepsilon/2, \lambda + \varepsilon/2)}(A|_{X_n}) i_{X_n}^*$. 
Then,
$$
\mbox{\rm Rank}(\cP_n) \geq q,
$$
and there exists $C\in \R_+$ such that, for $n$ large enough, 
\begin{equation}\label{eq:proj1exact}
 \|(\cP-\cP_n)\cP\|_{\cL(\cH, \QA)} \leq C \left\| \left( 1 - \Pi_{X_n}^{\QA}\right) \cP \right\|_{\cL(\cH, \QA)}.
\end{equation}
If we assume in addition that
$$  
\mbox{\rm (D2)} \quad \mbox{for $n$ large enough, $\mbox{\rm Rank}(\cP_n) = q$},
$$
then there exists $C \in \R_+$ such that, for $n$ large enough, 
\begin{equation}\label{eq:proj2exact}
 \|(\cP-\cP_n)\cP_n\|_{\cL(\cH,\QA)} \leq  C \left\| \left( 1 - \Pi_{X_n}^{\QA}\right) \cP \right\|_{\cL(\cH, \QA)},
\end{equation}
\begin{equation}\label{eq:eigenvalueexact}
 \dps \mathop{\max}_{\lambda_n \in \sigma(A_n)\cap (\lambda - \varepsilon/2, \lambda+\varepsilon/2)} |\lambda_n - \lambda| \leq C \left( \left\| \left( 1 - \Pi_{X_n}^{\QA} \right) \cP \right\|_{\cL(\cH, \QA)}  \right)^2. 
\end{equation}
\end{corollary}

\medskip

The estimates (\ref{eq:proj1exact}), (\ref{eq:proj2exact}) and (\ref{eq:eigenvalueexact}) are optimal. They are similar to the ones proved in \cite{Rappaz2,Mills2, Mills1}, but our assumptions on the sequence of discretized operators $A|_{X_n}$ are different. In \cite{Rappaz2}, these estimates are proved under the condition 
\begin{equation}\label{eq:Rappazcond}
 \delta(A, A|_{X_n}) \mathop{\longrightarrow}_{n\to\infty} 0,
\end{equation}
where
$$
 \delta(A, A|_{X_n}) := \mathop{\sup}_{\phi\in \DA, \; \|\phi\|_{\cH} + \|A\phi\|_{\cH} = 1} \mathop{\inf}_{\phi_n\in X_n} \|\phi - \phi_n\|_{\cH} + \|A\phi - A|_{X_n}\phi_n \|_{\cH}.
$$
In \cite{Mills2}, the assumptions are that $A$ is invertible and 
\begin{equation}\label{eq:Millscond}
 \mathop{\sup}_{v_n\in X_n} \mathop{\inf}_{w_n\in X_n}\frac{\left\| A^{-1}v_n -w_n\right\|_\QA}{\|v_n\|_{\QA}} \mathop{\longrightarrow}_{n\to\infty} 0.
\end{equation}
Each of the conditions (\ref{eq:Rappazcond}) and (\ref{eq:Millscond}) ensures that (D1) and (D2) hold for any discrete eigenvalue of~$A$. 
In the case when $A$ is semibounded, (C3) is automatically satisfied, so that our assumptions boil down to (A1), (D1) and (D2). These three conditions are weaker than those in~\cite{Rappaz2,Mills2, Mills1}, and more easy to check in some settings, as will be seen in Section~\ref{sec:supercell} on the example of the supercell method. 
On the other hand, when $A$ is not semibounded, the precise relationship between condition (C3) and (\ref{eq:Rappazcond}) and (\ref{eq:Millscond}) is still unclear to us. 

\section{Proof of Theorem~\ref{th:thmain}}\label{sec:proofmain}

\subsection{Proof of (\ref{eq:atteinte})}\label{sec:proof1}

Let us argue by contradiction and assume that there exists a subsequence $(\cT_{n_k})_{k\in\N}$ and $\eta>0$ such that $(\lambda-\eta, \lambda+\eta)\cap \sigma(A) = \{\lambda\}$ and 
\begin{equation}\label{eq:absurde}
\forall k\in \N, \; \mbox{\rm dist}\left( \lambda , \sigma(A_{n_k})\right) \geq \eta.
\end{equation}
Let $\psi\in\DA$ be a $\cH$-normalized eigenvector of $A$ associated with the discrete eigenvalue $\lambda$ and $\mu := \lambda + \eta/2$. As 
$(\mu-\frac{\eta}{2}, \mu + \frac{\eta}{2})\cap \sigma(A) = \emptyset$, it holds
$$
\alpha := \mathop{\min}_{\nu \in \sigma(A)} \frac{|\nu - \mu|}{1+|\nu|} >0.
$$
Let us consider the auxiliary problem
\begin{equation}\label{eq:pbcont1}
\left\{
\begin{array}{l}
 \mbox{find }u\in \QA\mbox{ such that}\\
\forall v\in \QA, \; (a-\mu m)(u,v) = (\lambda - \mu)m(\psi, v).
\end{array}
\right .
\end{equation}
The bilinear form $a-\mu m$ is continuous on $\QA$ and satisfies $\|a-\mu m\|_{\cL(\QA\times \QA)} \leq 1+|\mu|$. The linear form 
$f:\QA \ni v \mapsto (\lambda - \mu)m( \psi, v)$ 
is also continuous. Furthermore, as $\mu \notin \sigma(A)$, if $v\in \QA$ is such that $(a-\mu m)(v,w) = 0$ for all $w\in \QA$, then necessarily $v=0$. Lastly,
$$
\mathop{\inf}_{w\in\QA} \mathop{\sup}_{v\in\QA} \frac{|(a-\mu m)(v,w)|}{\|v\|_{\QA}\|w\|_{\QA}} \geq \mathop{\min}_{\nu \in \sigma(A)} \frac{|\nu - \mu|}{1 + |\nu|} = \alpha.
$$
Thus, applying Banach-Ne\v{c}as-Babu\v{s}ka's theorem (see Section~\ref{sec:appendix}), problem (\ref{eq:pbcont1}) is well-posed. Clearly, its unique solution is $u=\psi$.

Let us now introduce the following sequence of discretized problems for $k\in\N$:
\begin{equation}\label{eq:pbdisc1}
\left\{
\begin{array}{l}
\mbox{find }u_{n_k}\in X_{n_k}\mbox{ such that}\\
\forall v_{n_k}\in X_{n_k}, \; (a_{n_k} - \mu m_{n_k})(u_{n_k}, v_{n_k}) = (\lambda - \mu)m_{n_k}\left( \Pi_{X_{n_k}}^{\QA}\psi, v_{n_k}\right).
\end{array}
\right .
\end{equation}
From (\ref{eq:absurde}) and assumption (A3) (for $K=\left\{\mu\right\}$ and $\alpha_K=\eta/2$), we deduce the discrete inf-sup condition
$$
\forall k\in\N, \quad \mathop{\inf}_{w_{n_k}\in X_{n_k}} \mathop{\sup}_{v_{n_k}\in X_{n_k}} \frac{|(a_{n_k} - \mu m_{n_k})(v_{n_k}, w_{n_k})|}{\|v_{n _k}\|_{\QA} \|w_{n_k}\|_{\QA}} \geq c > 0.
$$
Thus, by Strang's lemma (see Section~\ref{sec:appendix}) and assumptions (A2), (A3) and (A4), for all 
$k\in\N$, 
\begin{eqnarray*} \|\psi - u_{n_k}\|_{\QA} 
& \leq & \frac{\eta}{2c}\mathop{\sup}_{v_{n_k}\in X_{n_k}} \frac{|m(\psi, v_{n_k}) - m_{n_k}(\Pi_{X_{n_k}}^\QA\psi, v_{n_k})|}{\|v_{n_k}\|_{\QA}} \\
&+& \inf_{w_{n_k}\in X_{n_k}}\left( \frac{c+1+|\mu|}{c}\|\psi -w_{n_k}\|_{\QA} + \frac{1}{c}\sup_{v_{n_k}\in X_{n_k}} \frac{|[(a_{n_k}-a) +\mu(m- m_{n_k})](w_{n_k},v_{n_k})|}{\|v_{n_k}\|_{\QA}}\right)\\
&\leq& \frac{\eta}{2c} \left( \left\|\psi - \Pi_{X_{n_k}}^\QA\psi \right\|_{\cH} + \kappa r_{n_k}^m\left(\Pi_{X_{n_k}}^{\QA}\psi\right) + s_{n_k}^m(\psi)\right) + \frac{c+1+|\mu|}{c}\left\|\psi -\Pi^\QA_{X_{n_k}} \psi\right\|_{\QA}\\
&& + \frac{1}{c}\left(\kappa r_{n_k}^a\left( \Pi_{X_{n_k}}^{\QA}\psi\right) + s_{n_k}^a(\psi) + |\mu|\kappa r_{n_k}^m\left(\Pi_{X_{n_k}}^{\QA} \psi\right) + |\mu|s_{n_k}^m(\psi)\right) \\
& \leq& C\left(  \left\|\psi -\Pi^\QA_{X_{n_k}} \psi\right\|_{\QA} + r_{n_k}^a(\psi) + r_{n_k}^m(\psi) + s_{n_k}^a(\psi) + s_{n_k}^m(\psi) \right),
\end{eqnarray*}
where $C \in \R_+$ is a constant independent of $k$. The above inequality implies that the sequence $(u_{n_k})_{k\in\N}$ strongly converges to $\psi$ in $\QA$, from which we infer that 
\begin{equation}\label{eq:convPIu}
\lim_{k \to \infty} \|\Pi_{X_{n_k}}^{\QA}\psi - u_{n_k}\|_{\QA} = 0.
\end{equation}
On the other hand, (\ref{eq:pbdisc1}) yields
$$
\forall v_{n_k} \in X_{n_k}, \; (a_{n_k} - \lambda m_{n_k})(u_{n_k}, v_{n_k})  = (\lambda - \mu)m_{n_k}\left( \Pi_{X_{n_k}}^{\QA}\psi - u_{n_k}, v_{n_k}\right).
$$
The above equality also reads 
$$
(A_{n_k} - \lambda)(\cM_{n_k}^{1/2}u_{n_k}) = (\lambda-\mu) {\cM}_{n_k}^{1/2}\left( \Pi_{X_{n_k}}^{\QA}\psi - u_{n_k} \right).
$$
It then follows from (A1), (A2) and (\ref{eq:convPIu}) that 
$$
\lim_{k \to \infty} \|(A_{n_k} - \lambda)(\cM_{n_k}^{1/2}u_{n_k})\|_{\cH} = 0 \qquad \mbox{ and } \qquad \liminf_{k \to \infty} \|  \cM_{n_k}^{1/2}u_{n_k} \|_\cH^2 \ge \gamma > 0, 
$$
which proves that $\dps \mbox{\rm dist}( \lambda,\sigma(A_{n_k})) \mathop{\longrightarrow}_{k\to\infty} 0$ and contradicts (\ref{eq:absurde}).

\subsection{Proof of  (\ref{eq:dim}) and (\ref{eq:proj1})}

By assumption (B1), the approximation $(\cT_n)_{n\in\N}$ is such that $\dps \mathop{\overline{\lim}}_{n\to\infty} \sigma(A_n) \cap (\lambda - \varepsilon, \lambda + \varepsilon)  = \{\lambda\}$. Hence, for $n$ large enough, 
$$
\sigma(A_n) \cap ( (\lambda - 2\varepsilon/3, \lambda - \varepsilon/3) \cup (\lambda + \varepsilon/3, \lambda +2\varepsilon/3)) = \emptyset,
$$
so that the circle $\cC$ in the complex plane centered at $\lambda$ and of radius $\varepsilon/2$ is such that $\mbox{\rm dist}(\cC, \sigma(A_n)) \geq \varepsilon/6$. 
This implies in particular that, for $n$ large enough, 
$$
\cP_n = \frac{1}{2i\pi} \oint_{\cC} i_{X_n}\cM_n^{-1/2}(z-A_n)^{-1}\cM_n^{1/2}i_{X_n}^*\,dz.
$$
Consequently, for all $\psi\in \mbox{\rm Ran}(\cP)$, it holds
$$
(\cP - \cP_n)\psi = \frac{1}{2i\pi} \oint_\cC \left( (z-A)^{-1}\psi - i_{X_n}\cM_n^{-1/2}(z-A_n)^{-1}\cM_n^{1/2}i_{X_n}^*\psi\right)\,dz.
$$
In the following, $C$ will denote a constant independent of $n\in\N^*$ and $z\in\cC$, which may change along the calculations.

For $z \in \cC$, we introduce the auxiliary problem
$$
\left\{
\begin{array}{l}
 \mbox{find }u^z\in\QA\mbox{ such that }\\
\forall v\in \QA, \; (z m -a)(u^z, v) = m\left( \psi, v \right),\\
\end{array}
 \right .
$$
whose unique solution is $u^z = (z-A)^{-1}\psi =  \frac{\psi}{z-\lambda}$, since $\psi \in \mbox{\rm Ran}(\cP)$. We also introduce the discretized problem
$$
\left\{
\begin{array}{l}
 \mbox{find }u^z_n\in X_n\mbox{ such that}\\
\forall v_n\in X_n, \; (zm_n - a_n)(u_n^z, v_n) = m_n \left( \Pi_{X_n}^{\cH}\psi, v_n \right),\\
\end{array}
\right.
$$
whose unique solution is $u_n^z = i_{X_n} \cM_n^{-1/2}(z-A_n)^{-1} \cM_n^{1/2} i_{X_n}^*\psi$. From assumption (A3), since $\cC$ is a compact subset of $\C$, 
there exists $c>0$ such that for all $z\in \cC$ and $n\in\N$,
$$
\mathop{\inf}_{w_n\in X_n} \mathop{\sup}_{v_n\in X_n} \frac{|(a_n - zm_n)(v_n,w_n)|}{\|v_n\|_{\QA} \|w_n\|_{\QA}} \geq c.
$$
Reasoning as in Section~\ref{sec:proof1}, we infer from Strang's lemma, assumptions (A2)-(A4) 
and the fact that $r_n^a, r_n^m, s_n^a$ and $s_n^m$ are semi-norms, that for all $z \in \cC$,
\begin{eqnarray*}
 \|u^z - u_n^z\|_{\QA} & \leq & \dps \frac{1}{c} \mathop{\sup}_{v_n\in X_n} \frac{|m(\psi, v_n) - m_n(\Pi_{X_n}^{\cH}\psi, v_n)|}{\|v_n\|_{\QA}}\\
&& + \mathop{\inf}_{w_n \in X_n} \left( \frac{c+1+|z|}{c}\|u^z - w_n\|_{\QA} + \frac{1}{c} \mathop{\sup}_{v_n\in X_n} \frac{|[(a_n-a) + z (m_n-m)](w_n, v_n)|}{\|v_n\|_{\QA}}\right)\\
& \leq & \frac{\kappa}{c} r_n^m(\psi) + \frac{\Gamma}{c}\left\|\psi - \Pi_{X_n}^{\cH}\psi\right\|_{\cH} 
+ \frac 1c s^m_n(\psi) + \frac \Gamma c \left\|\Pi_{X_n}^{\QA}\psi - \Pi_{X_n}^{\cH}\psi\right\|_{\cH} \\
&&  + \frac{c+1+|z|}{c}\left\| u^z - \Pi_{X_n}^{\QA}u^z \right\|_{\QA} \\
&&  + \frac{1}{c}\left(\kappa r_n^a\left(\Pi_{X_n}^{\QA}u^z\right) + s_n^a(u^z) + |z| \kappa r_n^m\left(\Pi_{X_n}^{\QA}u^z\right) + |z| s_n^m(u^z) \right)\\
& \leq &  \frac{\kappa}{c} r_n^m(\psi)   + \frac{3\Gamma}{c}\left\|\psi - \Pi_{X_n}^{\QA}\psi\right\|_{\QA} + \frac{c+(1+|z|)(1+\kappa^2)}{c}\left\| u^z - \Pi_{X_n}^{\QA}u^z \right\|_{\QA} \\
&&  + \frac{1}{c}\left(\kappa r_n^a(u^z) + s_n^a(u^z) + |z| \kappa r_n^m(u^z) + |z| s_n^m(u^z)\right) \\
& \le  & C \left( \left\| \left( 1 - \Pi_{X_n}^\QA\right) \psi \right\|_{\QA} + r_n^a(\psi) + r_n^m(\psi) + s_n^a(\psi) + s_n^m(\psi) \right),
\end{eqnarray*}
since $u^z= \frac{\psi}{z-\lambda}$. Thus, for all $z\in\cC$, 
\begin{eqnarray*}
 && \!\!\!\!\!\!\!\!\! \left\| (z-A)^{-1}\psi - i_{X_n} \cM_n^{-1/2} (z-A_n)^{-1} \cM_n^{1/2} i_{X_n}^* \psi \right\|_{\QA} \\
&& \qquad \qquad \leq C \left( \left\| \left( 1 - \Pi_{X_n}^{\QA}\right) \psi \right\|_{\QA} + r_n^a(\psi) + r_n^m(\psi) + s_n^a(\psi) + s_n^m(\psi) \right).\\
\end{eqnarray*}
Since $\cC$ is of finite length, we obtain that, for $n$ large enough, for all $\psi \in \mbox{\rm Ran}(\cP)$,
$$
\|(\cP - \cP_n)\psi \|_{\QA} \leq C \left( \left\| \left( 1 - \Pi_{X_n}^{\QA}\right)\psi \right\|_{\QA} + r_n^a(\psi) + r_n^m(\psi) + s_n^a(\psi) + s_n^m(\psi)\right),
$$
which readily leads to (\ref{eq:proj1}). 

Let us finally consider a $\cH$-orthonormal basis $(\zeta_1, \cdots , \zeta_q)$ of $\mbox{\rm Ran}(\cP)=\mbox{Ker}(\lambda-A)$. Since for all $1\leq i \leq q$, $\dps \cP_n \zeta_i \mathop{\longrightarrow}_{n\to\infty} \cP\zeta_i = \zeta_i$ strongly in $\cH$, the family $(\cP_n \zeta_1,\cdots, \cP_n \zeta_q)$ is free for $n$ large enough, so that $\mbox{Rank}(\cP_n) \geq q$. 

\subsection{Proof of (\ref{eq:proj2}) and (\ref{eq:eigenvalue})}

We just have shown that for all $1\leq i \leq q$, $\dps \cP_n \zeta_i \mathop{\longrightarrow}_{n\to\infty} \cP\zeta_i = \zeta_i$ strongly in $\cH$. 
Under the additional assumption that, for $n$ large enough, $\mbox{\rm Rank}(\cP_n)=q$, this implies that there exists
 $n_0 \in \N$, such that, for $n \ge n_0$,  $(\cP_n\zeta_1, \cdots, \cP_n \zeta_q)$ forms a basis of $\mbox{\rm Ran}(\cP_n)$, 
with
$$
\mathop{\min}_{1\leq i \leq q} \|\cP_n \zeta_i\|^2_{\cH} \geq \frac{3}{4} \quad \mbox{\rm and} \quad \mathop{\max}_{1\leq i,j \leq q,\; i\neq j} \left| m\left( \cP_n\zeta_i, \zeta_j\right)\right| \leq \frac{1}{4q}.
$$
Thus, any $\xi_n \in \mbox{\rm Ran}(\cP_n)$ can be decomposed as
$$
\xi_n = \sum_{i=1}^q \alpha^i(\xi_n) \cP_n \zeta_i,
$$
the coefficients $(\alpha^1(\xi_n), \cdots, \alpha^q(\xi_n))$ of $\xi_n$ in the basis $(\cP_n \zeta_1, \cdots, \cP_n \zeta_q)$ being such that
$$
\mathop{\max}_{1\leq i \leq q} |\alpha^i(\xi_n)| \leq 2 \|\xi_n\|_{\cal H}. 
$$  
We have
\begin{eqnarray*}
 \cP\xi_n-\xi_n &= & \sum_{i=1}^q \alpha^i(\xi_n)\left( \sum_{j=1}^q m\left( \cP_n \zeta_i, \zeta_j\right) \zeta_j - \cP_n \zeta_i\right) \\
& = & \sum_{i=1}^q \alpha^i(\xi_n) \left( \sum_{j\neq i} m\left( \cP_n \zeta_i - \zeta_i, \zeta_j\right) \zeta_j - (\cP_n \zeta_i- \zeta_i) + m\left( \cP_n \zeta_i - \zeta_i, \zeta_i\right) \zeta_i \right),
\end{eqnarray*}
and we deduce from (\ref{eq:proj1}) that for all $1\leq i \leq q$,
$$
\|\zeta_i - \cP_n\zeta_i\|_\QA \leq C \left( \left\| \left( 1 - \Pi_{X_n}^\QA\right)\cP\right\|_{\cL(\cH, \QA)} + \cR_n^a + \cR_n^m + \cS_n^a + \cS_n^m\right).
$$
Hence, 
$$
\forall \xi_n \in \mbox{\rm Ran}(\cP_n), \quad \|\cP\xi_n-\xi_n\|_{\QA} \leq C\left( \left\| \left( 1 - \Pi_{X_n}^\QA\right) \cP \right\|_{\cL(\cH, \QA)} + \cR_n^a + \cR_n^m + \cS_n^a + \cS_n^m\right) \|\xi_n\|_{\cal H}, 
$$
where the constant $C$ is independent of $n$. Besides, it also follows from (A2) and the definition of $\cP_n$ that
$$
\forall n \in \N, \quad \|\cP_n\|_{{\cal L}(\cH)} \le \sqrt{\frac\Gamma\gamma}.
$$
Therefore,
\begin{eqnarray*}
\|(\cP-\cP_n)\cP_n\|_{\cL(\cH,\QA)} &\leq&  \sup_{\xi_n \in {\rm Ran}(\cP_n)\setminus \left\{0\right\}}  \frac{\|\cP\xi_n-\xi_n\|_{\QA}}{ \|\xi_n\|_{\cH}} \;  \|\cP_n\|_{{\cal L}(\cH)} \\
&\leq& C\left( \left\| \left( 1 - \Pi_{X_n}^\QA\right) \cP \right\|_{\cL(\cH, \QA)} + \cR_n^a + \cR_n^m + \cS_n^a + \cS_n^m\right) ,
\end{eqnarray*}
and (\ref{eq:proj2}) is proved.

\medskip

For each $n$ large enough, let $(\psi_n,\lambda_n) \in X_n \times \R$ be a solution to the generalized eigenvalue problem (\ref{eq:geneig}) such that $\lambda_n \in (\lambda - \varepsilon/2, \lambda + \varepsilon/2)$, ${\phi}_n = \frac{\psi_n}{\|\psi_n\|_\cH}$, and $\chi_n = \frac{\cP \psi_n}{\|\cP\psi_n\|_{\cH}} = \frac{\cP \phi_n}{\|\cP\phi_n\|_{\cH}}$. It follows from (\ref{eq:proj2}) that
$$
\| \cP\phi_n-\phi_n\|_\cH \le \|\cP\phi_n-\phi_n\|_\QA \le C\left( \left\| \left( 1 - \Pi_{X_n}^\QA\right) \cP \right\|_{\cL(\cH, \QA)} + \cR_n^a + \cR_n^m + \cS_n^a + \cS_n^m\right) \mathop{\longrightarrow}_{n\to\infty} 0,
$$
from which we infer that $\|\cP\phi_n\|_\cH \to 1$, $(\phi_n)_{n \in \N}$ is bounded in $\QA$, $\|\phi_n-\chi_n\|_\QA \to 0$, and
\begin{eqnarray*}
\|\chi_n-\phi_n\|_\QA &\le& \left\| \frac{\cP\phi_n}{\|\cP\phi_n\|_\cH}-\cP\phi_n \right\|_\QA + \|\cP\phi_n-\phi_n\|_\QA  \nonumber \\
&\le& \|\cP\phi_n\|_\cH^{-1} \|\phi_n-\cP\phi_n\|_\cH \|\cP\phi_n\|_\QA +
\|\cP\phi_n-\phi_n\|_\QA  \nonumber\\ 
&\le&  C\left( \left\| \left( 1 - \Pi_{X_n}^\QA\right) \cP \right\|_{\cL(\cH, \QA)} + \cR_n^a + \cR_n^m + \cS_n^a + \cS_n^m\right). 
\end{eqnarray*}
Besides, it holds
\begin{eqnarray*}
|\lambda_n-\lambda| &=& |a_n(\psi_n,\psi_n)-a(\chi_n,\chi_n)| \nonumber \\
&\leq& \left|\frac{a_n(\phi_n,\phi_n)}{m_n(\phi_n,\phi_n)}-a(\phi_n,\phi_n)\right|+|a(\phi_n,\phi_n)-a(\chi_n,\chi_n)|.
\end{eqnarray*}
On the one hand, we have
\begin{eqnarray*}
 |a(\phi_n, \phi_n) - a(\chi_n, \chi_n)| & = & |a(\phi_n - \chi_n, \phi_n - \chi_n) + 2a(\chi_n, \phi_n - \chi_n)|\nonumber \\
& = & |a(\phi_n - \chi_n, \phi_n - \chi_n) + 2\lambda m( \chi_n, \phi_n - \chi_n)|\nonumber \\
& = & |a(\phi_n - \chi_n, \phi_n - \chi_n) - \lambda \|\chi_n - \phi_n\|_{\cH}^2|\nonumber \\
& \le &  C\|\phi_n - \chi_n\|_{\QA}^2.
\end{eqnarray*}
On the other hand,
\begin{eqnarray*}
|(a-a_n)(\phi_n, \phi_n)| & \leq  & |(a-\widetilde{a}_n)(\phi_n, \phi_n)| + |(\widetilde{a}_n-a_n)(\phi_n, \phi_n)|\\
& \leq & r_n^a(\phi_n)^2 + \left|(\widetilde{a}_n-a_n)\left(\phi_n - \Pi_{X_n}^{\QA} \chi_n, \phi_n - \Pi_{X_n}^{\QA}\chi_n\right)\right|\\
&& + 2 \left|(\widetilde{a}_n-a_n)\left(\Pi_{X_n}^{\QA} \chi_n, \phi_n - \Pi_{X_n}^{\QA}\chi_n\right)\right|\\
&& + \left|(\widetilde{a}_n-a_n)\left(\Pi_{X_n}^{\QA} \chi_n, \Pi_{X_n}^{\QA}\chi_n\right)\right|\\
& \leq & \left( r_n^a(\chi_n) + \kappa \|\phi_n - \chi_n\|_{\QA}\right)^2 + (\Gamma+\kappa^2) \left\| \phi_n - \Pi_{X_n}^\QA \chi_n \right\|_{\QA}^2 \\
&& + \left( 2 \left\|\phi_n - \Pi_{X_n}^\QA \chi_n \right\|_{\QA} +\left\|\Pi_{X_n}^\QA \chi_n \right\|_{\QA} \right) s_n^a(\chi_n)\\
& \leq & C \left[ \left(   r_n^a(\chi_n) + \|\phi_n - \chi_n\|_{\QA} + \left\| \left( 1 - \Pi_{X_n}^{\QA} \right) \chi_n \right\|_{\QA} \right)^2 + s_n^a(\chi_n) \right] \\
& \leq &  C \left[\left( \left\| \left( 1 - \Pi_{X_n}^\QA\right) \cP \right\|_{\cL(\cH, \QA)} + \cR_n^a + \cR_n^m + \cS_n^a + \cS_n^m\right)^2 + \cS_n^a  \right],
\end{eqnarray*}
and a similar calculation leads to 
\begin{eqnarray*}
|m_n(\phi_n, \phi_n) -1| & = & |m_n(\phi_n, \phi_n) -m(\phi_n, \phi_n)|,\\
& \leq &  C \left[\left( \left\| \left( 1 - \Pi_{X_n}^\QA\right) \cP \right\|_{\cL(\cH, \QA)} + \cR_n^a + \cR_n^m + \cS_n^a + \cS_n^m\right)^2 + \cS_n^m  \right].
\end{eqnarray*}
Consequently,
\begin{eqnarray*}
 \left| \frac{a_n(\phi_n,\phi_n)}{m_n(\phi_n, \phi_n)} - a(\phi_n, \phi_n) \right|  &\leq& \frac{|(a-a_n)(\phi_n, \phi_n)|}{m_n(\phi_n, \phi_n)} + |a(\phi_n, \phi_n)| \left| \frac{m_n(\phi_n, \phi_n) -1}{m_n(\phi_n, \phi_n)}\right| \nonumber \\
& \leq & \gamma^{-1} \left( |(a-a_n)(\phi_n, \phi_n)| +  |a(\phi_n, \phi_n)| \, |m_n(\phi_n, \phi_n) -1| \right) \nonumber \\
& \leq &  C \left[\left( \left\| \left( 1 - \Pi_{X_n}^\QA\right) \cP \right\|_{\cL(\cH, \QA)} + \cR_n^a + \cR_n^m + \cS_n^a + \cS_n^m\right)^2 + \cS_n^a  + \cS_n^m  \right]. \nonumber
\end{eqnarray*}
Collecting the above results, we obtain
$$
|\lambda - \lambda_n| \leq C \left[\left( \left\| \left( 1 - \Pi_{X_n}^\QA \right) \cP \right\|_{\cL(\cH, \QA)} + \cR_n^a + \cR_n^m\right)^2 + \cS_n^a + \cS_n^m \right],
$$
which proves estimate (\ref{eq:eigenvalue}).

\section{Application to the supercell method}\label{sec:supercell}

The aim of this section is to show that the theoretical framework presented in Section~\ref{sec:theory} can be applied to the numerical analysis of the 
supercell method for perturbed periodic Schr\"odinger operators. 

\medskip

Note that the supercell method was previously studied from a mathematical viewpoint by Soussi~\cite{Soussi}, for the special case of a 
two-dimensional periodic Schr\"odinger operator in the presence of a compactly supported perturbation $W$ of the form $W(x) = w\un_{\Omega}(x)$, where 
$w$ is a real constant and $\Omega$ a bounded domain of $\R^2$.

\subsection{The supercell method with exact integration}\label{sec:supercelldescription}

Let $\cR$ be a periodic lattice of $\R^d$, $\cR^*$ its reciprocal lattice and $\Gamma$ a unit cell of $\cR$ such that $0$ is in the interior of $\Gamma$. Typically, in the case of the cubic lattice $\cR = \Z^d$, $\cR^* = 2\pi\Z^d$ and $\Gamma = (-1/2, 1/2]^d$ is an admissible unit cell.  

Let us introduce the perturbed periodic Schr\"odinger operator 
$$
A:= -\Delta + V_{\rm per} +W,
$$
where $\Delta$ is the Laplace operator, $V_{\rm per}$ a real-valued $\cR$-periodic function of $L^p_{\rm loc}(\R^d)$, with $p=2$ if $d\leq 3$, $p>2$ if $d=4$
 and $p=d/2$ for $d\geq 5$, and $W\in L^{\infty}(\R^d)$ a real-valued function such that $\dps W(x) \mathop{\longrightarrow}_{|x|\to\infty} 0$. 

The operator $A$ is self-adjoint and bounded from below on $\cH:= L^2(\R^d)$, 
endowed with its natural inner product
$$
\forall \phi, \psi \in \cH, \quad m(\phi,\psi) := \int_{\R^d} \phi \psi,
$$
with domain $D(A) = H^2(\R^d)$ and form domain $\QA = H^1(\R^d)$. The associated bilinear form $a(\cdot, \cdot)$ is defined by
$$
\forall \phi,\psi \in \QA, \quad a(\phi, \psi) := \int_{\R^d} \nabla \phi \cdot \nabla \psi + \int_{\R^d} (V_{\rm per} + W) \phi \psi.
$$
We denote by $A^0:= -\Delta + V_{\rm per}$ the corresponding periodic Schr\"odinger operator on $L^2(\R^d)$. 

\medskip

The supercell method is the current state-of-the-art technique in solid state physics to compute the spectrum of the operator $A$.
For $L\in\N^*$, we denote by $\Gamma_L:= L\Gamma$ the supercell of size $L$ and set
\begin{eqnarray*}
L^2_{\rm per}(\Gamma_L) & := & \left\{ u_L \in L^2_{\rm loc}(\R^d) \; | \; u_L \; L\cR\mbox{-periodic} \right\},\\
H^1_{\rm per}(\Gamma_L) & := & \left\{ u_L \in L^2_{\rm per}(\Gamma_L) \; | \; \nabla u_L \in \left( L^2_{\rm per}(\Gamma_L)\right)^d \right\},\\
C^0_{\rm per}(\Gamma_L) & := & \left\{  u_L \in C^0(\R^d) \; |\; u_L \; L\cR\mbox{-periodic} \right\},\\
L^{\infty}_{\rm per}(\Gamma_L) & := & \left\{ u_L \in L^{\infty}(\R^d)\; |\; u_L \; L\cR\mbox{-periodic} \right\}.
\end{eqnarray*}
For $u_L \in L^2_{\rm per}(\Gamma_L)$ and $k\in L^{-1}\cR^*$, we denote by
$$
\uh_L(k):= \frac{1}{|\Gamma_L|^{1/2}}\int_{\Gamma_L} u_L(x)e^{ik\cdot x}\, dx
$$
the Fourier coefficient of $u_L$ corresponding to the $k$ mode. For $r \in \R$, the Sobolev space $H^r_{\rm per}(\Gamma_L)$ can be defined as 
$$
H^r_{\rm per}(\Gamma_L):= \left\{ u_L \in L^2_{\rm per}(\Gamma_L)\; | \; \sum_{k\in L^{-1}\cR^*} \left( 1 + |k|^2\right)^{r} \left| \uh_L(k)\right|^2 < \infty \right\}.
$$
The supercell method relies on the resolution of the following (non-consistent and non-conforming) eigenvalue problem:
$$
\left\{
\begin{array}{l}
 \mbox{find }(u_{L,N}, \lambda_{L,N}) \in Y_{L,N} \times \R \mbox{ such that}\\
\forall v_{L,N} \in Y_{L,N}, \; \ah_L(u_{L,N}, v_{L,N}) = \lambda_{L,N} \mh_L\left( u_{L,N}, v_{L,N}\right),
\end{array}
\right .
$$
where
\begin{eqnarray*}
\forall u_L,v_L  \in L^2_{\rm per}(\Gamma_L), &\;&  \mh_L\left( u_L, v_L \right)  := \int_{\Gamma_L}u_Lv_L,\\
\forall u_L,v_L \in H^1_{\rm per}(\Gamma_L), &\;& \ah_L(u_L,v_L)  :=  \int_{\Gamma_L} \nabla u_L \cdot \nabla v_L + \int_{\Gamma_L}(V_{\rm per} + W)u_Lv_L,
\end{eqnarray*}
and $Y_{L,N}$ is a finite dimensional subspace of $H^1_{\rm per}(\Gamma_L)$. 

We set $H_{L,N} = H_L|_{Y_{L,N}}$, where $H_L$ denotes the unique self-adjoint operator on $L^2_{\rm per}(\Gamma_L)$ associated with the quadratic form $\ah_L$. We have $D(H_L) = H^2_{\rm per}(\Gamma_L)$ and
$$
\forall u_L \in H^2_{\rm per}(\Gamma_L), \quad H_L u_L = -\Delta u_L + (V_{\rm per} +W_L)u_L,
$$
where $W_L\in L^{\infty}_{\rm per}(\Gamma_L)$ denotes the $L\cR$-periodic extension of $W|_{\Gamma_L}$.
\medskip

For the sake of clarity, our analysis will be restricted to the case of the cubic lattice $\cR = \Z^d$ and the planewave discretization method, for which
$$
Y_{L,N} := \left\{ \sum_{k\in L^{-1} \cR^\ast \;|\; |k|\leq 2\pi N L^{-1}} c_k e_{L,k} \; | \; \forall k,\; c_{-k} = c_k^* \right\},
$$
where $e_{L,k}(x) := |\Gamma_L|^{-1/2}e^{ik\cdot x}$. We denote by $\Pi_{L,N}$ the orthogonal projection of $L^2_{\rm per}(\Gamma_L)$ on $Y_{N,L}$ for the $L^2_{\rm per}(\Gamma_L)$ inner product (actually $\Pi_{L,N}$ is also the orthogonal projection of $H^s_{\rm per}(\Gamma_L)$ on $Y_{N,L}$ for the $H^s_{\rm per}(\Gamma_L)$ inner product, for any $s \in \R$). 

The discretization spaces $Y_{L,N}$ possess the following properties:
$$
\forall u_{L,N} \in Y_{L,N},\; \Pi_{L,N}\left( -\Delta u_{L,N} \right) = -\Delta u_{L,N},
$$
and for all real numbers $r$ and $s$ such that $0\leq r \leq s$, there exists a constant $C \in \R_+$ such that for all $L\in\N^*$ and 
all $u_L\in H^s_{\rm per}(\Gamma_L)$, 
\begin{equation}\label{eq:Jackson}
\|u_L - \Pi_{L,N}u_L\|_{H^r_{\rm per}(\Gamma_L)} \leq C\left( \frac{L}{N}\right)^{s-r} \|u_L\|_{H^s_{\rm per}(\Gamma_L)}.
\end{equation}

As in \cite{CEM}, we will assume that $V_{\rm per}$ belongs to the functional space $\cZ_{\rm per}(\Gamma)$ (denoted by $\cM_{\rm per}(\Gamma)$ in \cite{CEM}), defined by
$$
\cZ_{\rm per}(\Gamma) := \left\{ V\in L^2_{\rm per}(\Gamma)\; |\; \|V\|_{\cZ_{\rm per}(\Gamma)} := \mathop{\sup}_{L\in\N^*} \mathop{\sup}_{w_L\in H^1_{\rm per}(\Gamma_L)\setminus \{0\}} \frac{\|Vw_L\|_{L^2_{\rm per}(\Gamma_L)}}{\|w_L\|_{H^1_{\rm per}(\Gamma_L)}} < +\infty \right\}.
$$
The space $\cZ_{\rm per}(\Gamma)$ is a normed space and the space of the $\cR$-periodic functions of class $C^{\infty}$ is dense in $\cZ_{\rm per}(\Gamma)$. 

\medskip

Our main result concerning the supercell method in the absence of numerical integration error is the following:
\begin{theorem}\label{th:supercell}
Assume that $V_{\rm per} \in \cZ_{\rm per}(\Gamma)$ and that $W \in L^\infty(\R^d)$ with $\dps W(x) \mathop{\longrightarrow}_{|x| \to \infty} 0$. Let $(N_L)_{L\in\N^*}$ be a sequence of integers such that $\dps \frac{N_L}{L} \mathop{\longrightarrow}_{L\to\infty} + \infty$. Then, 
\begin{itemize}
 \item[1.] \bfseries Absence of pollution \normalfont \itshape
\end{itemize}
\begin{equation}\label{eq:nopollutionsc}
\mathop{\lim}_{L\to\infty} \sigma(H_{L,N_L}) = \sigma(A).
\end{equation}
\begin{itemize}
 \item[2.] \bfseries A priori error estimates 
\end{itemize}
Assume that, in addition, $V_{\rm per} \in H^{r-2}_{\rm per}(\Gamma)$ and $W \in H^{r-2}(\R^d)$, for some $r \ge 2$. Let $\lambda$ be a discrete eigenvalue of $A$ and $\varepsilon >0$ be such that $\sigma(A) \cap (\lambda - \varepsilon, \lambda +\varepsilon) = \{\lambda\}$. 
Let $\cP := \un_{\left\{\lambda\right\}}(A)$ be the $L^2(\R^d)$-orthogonal projection onto the eigenspace of $A$ associated with $\lambda$ and 
$\fP_L := \un_{(\lambda - \varepsilon/2, \lambda +\varepsilon/2)}(H_{L,N_L})$ the $L^2_{\rm per}(\Gamma_L)$-orthogonal spectral projection of $H_{L,N_L}$ associated with the eigenvalues belonging to the interval $(\lambda-\varepsilon/2, \lambda+\varepsilon/2)$. Consider finally a sequence of cut-off functions $(\chi_L)_{L \in \N^\ast}$ such that
\begin{equation}\label{eq:prop_chiL}
0\leq \chi_L\leq 1 \mbox{ on  } \R^d, \; \chi_L  = 1 \mbox{ on } \Gamma_L,
\; \mbox{\rm Supp}(\chi_L)\subset (L+\sqrt{L})\Gamma, \; \|\nabla \chi_L\|_{L^{\infty}} \leq c,
\end{equation}
for some constant $c\in\R_+$ independent of $L\in\N^*$. 

Then, $\mbox{\rm Ran}(\cP)\subset H^r(\R^d)$, and there exists $C,\delta>0$ such that for $L$ large enough,
\begin{equation}\label{eq:dimequalitysc}
\mbox{\rm Tr}(\cP) = \mbox{\rm Tr}(\fP_L),
\end{equation}
\begin{equation}\label{eq:proj1sc}
 \mathop{\sup}_{\psi\in {\rm Ran}(\cP), \; \|\psi\|_{L^2(\R^d)} = 1} \mathop{\inf}_{u_L \in {\rm Ran}(\fP_L)} \|\psi - \chi_Lu_L\|_{H^1(\R^d)} \leq C\left( e^{-\delta L} + \left(\frac{L}{N_L}\right)^{r-1}\right),
\end{equation}  
\begin{equation}\label{eq:proj2sc}
 \mathop{\sup}_{u_L\in {\rm Ran}(\fP_L), \; \|u_L\|_{L^2_{\rm per}(\Gamma_L)} = 1} \mathop{\inf}_{\psi \in {\rm Ran}(\cP)} \|\psi - \chi_Lu_L\|_{H^1(\R^d)} \leq C\left( e^{-\delta L} + \left(\frac{L}{N_L}\right)^{r-1}\right),
\end{equation}  
\begin{equation}\label{eq:eigenvaluesc}
\mathop{\max}_{\lambda_L \in \sigma(H_{L,N_L})\cap (\lambda-\varepsilon/2, \lambda+\varepsilon/2)} |\lambda_L - \lambda| \leq C\left(e^{-\delta L} + \left(\frac{L}{N_L}\right)^{r-1}\right)^2.
\end{equation}
\end{theorem}

\subsection{The supercell method with numerical integration}\label{sec:supercellnumint}

In general, the computation of the integral $\int_{\Gamma_L} (V_{\rm per} + W)u_Lv_L$ with $u_L, v_L \in Y_{L, N_L}$ cannot be carried out explicitly, and a numerical integration procedure is needed. We assume in this section that $V_{\rm per}$ and $W$ are continuous functions.

For $M\in\N^*$ and $u_L\in C^0_{\rm per}(\Gamma_L)$, we denote by $\uh_L^{FFT,M}$ the discrete Fourier transform of $u_L$ on the cartesian grid $\cG_{L,M}:= \frac{L}{M}\Z^d$. 
Recall that if 
$$
u_L = \sum_{k \in L^{-1}\cR^*} \uh_L(k) e_{L,k}, 
$$
the discrete Fourier transform of $u_L$ is the $ML^{-1}\cR^*$-periodic sequence $\uh_L^{FFT,M} = \left( \uh_L^{FFT,M}(k) \right)_{k\in L^{-1}\cR^*}$ where
$$
\uh_L^{FFT,M}(k) = \frac{1}{M^d} \sum_{x\in \cG_{L,M} \cap \Gamma_L} u_L(x) e^{-ik\cdot x} = |\Gamma_L|^{-1/2} \sum_{K\in L^{-1}\cR^*} \uh_L(k+MK).
$$
We now introduce the subspaces
$$
W_{L,M}^{1D}:= \left|
\begin{array}{ll}
\dps  \mbox{\rm Span}\left\{ e^{ily} \; |\; l\in 2\pi L^{-1}\Z, \; |l|\leq \frac{2\pi}{L}\left(\frac{M-1}{2}\right) \right\} & \quad (M \mbox{ odd}),\\
\dps \mbox{\rm Span}\left\{ e^{ily} \; |\; l\in 2\pi L^{-1}\Z, \; |l|\leq \frac{2\pi}{L}\left(\frac{M-1}{2}\right) \right\} \oplus \C\left( e^{i\pi My/L} + e^{-i\pi My/L}\right) & \quad (M \mbox{ even}),\\
\end{array}
\right .
$$
and denote by $W_{L,M}$ the $d$-tensor product space $W_{L,M}:= W_{L,M}^{1D} \otimes \cdots \otimes W_{L,M}^{1D}$. 
In particular, when $M$ is odd, 
$$
W_{L,M}= \mbox{\rm Span} \left\{ e_{L,k} , \; k \in L^{-1}\cR^*, \; |k|_{\infty} \leq 2\pi L^{-1} \left( \frac{M-1}{2}\right) \right\}.
$$
It is then possible to define the interpolation projector $\cI_{L,M}$ from $C^0_{\rm per}(\Gamma_L)$ onto $W_{L,M}$ by $\left[\cI_{L,M}(u_L)\right](x) = u_L(x)$ for all $x\in \cG_{L,M}$. 
In particular, when $M$ is odd, we have the simple relation
$$
\cI_{L,M}(u_L) = |\Gamma_L|^{1/2} \sum_{k\in L^{-1}\cR^* \; |\; |k|_{\infty} \leq 2\pi L^{-1} \left( \frac{M-1}{2}\right)} \uh_L^{FFT,M}(k) e_{L,k}.
$$
It is easy to check that if the function $u_L$ is real-valued, then so is the function $\cI_{L,M}(u_L)$. 

Besides, when $M\geq 4N+1$, it holds that for all $u_L, v_L \in Y_{L,N}$, 
$$
\int_{\Gamma_L} \cI_{L,M}(V_L u_L v_L) = \int_{\Gamma_L} \cI_{L,M}\left( V_L \right)u_L v_L,
$$
for any $V_L\in L^2_{\rm per}(\Gamma_L)$. 

\medskip

The supercell method with numerical integration then consists in considering the following eigenvalue problem for a given $M\geq 4N+1$,
$$
\left\{
\begin{array}{l}
 \mbox{find }(u_{L,N}, \lambda_{L,N}) \in Y_{L,N} \times \R \mbox{ such that}\\
\forall v_{L,N} \in Y_{L,N}, \; \ah_{L,M}(u_{L,N}, v_{L,N}) = \lambda_{L,N} \mh_L\left( u_{L,N}, v_{L,N}\right),\\
\end{array}
\right .
$$
where
\begin{eqnarray*}
\forall u_L,v_L \in H^1_{\rm per}(\Gamma_L), & \;& \ah_{L,M}(u_L,v_L) : =  \int_{\Gamma_L} \nabla u_L \cdot \nabla v_L + \int_{\Gamma_L}\cI_{L,M}(V_{\rm per} + \widetilde{W}_L)u_Lv_L,\\
\end{eqnarray*} 
and where $\widetilde{W}_L$ is the $L\cR$-periodic extension of $\xi_L W|_{\Gamma_L}$, $\xi_L$ being a $C^{[r-1]}(\R^d)$ cut-off function such that $0\leq \xi_L \leq 1$, $\xi_L = 1$ on $\Gamma_{L-1}$, $\mbox{\rm Supp}(\xi_L) \subset \left( L - 1/2\right) \Gamma$, and the sequences $\left( \|\partial^\alpha \xi_L\|_{L^\infty(\R^d)}\right)_{L\in\N^*}$ are uniformly bounded in $L$, for all $|\alpha| \le [r-1]$ (here and above, $[r-1]$ denotes the integer part of $r-1$).

\medskip

As in the preceding section, we denote by $H_{L,N,M} = H_{L,M}|_{Y_{L,N}}$, where $H_{L,M}$ is the unique self-adjoint operator on $L^2_{\rm per}(\Gamma_L)$ with domain $D(H_{L,M}) = H^2_{\rm per}(\Gamma_L)$ associated with the quadratic form $\ah_{L,M}$.

\begin{theorem}\label{th:supercellint}
Let $(N_L)_{L\in\N^*}$ and $(G_L)_{L\in\N^*}$ be sequences of integers such that 
$\dps \frac{N_L}{L} \mathop{\longrightarrow}_{L\to\infty} + \infty$ and $\dps G_L\mathop{\longrightarrow}_{L\to\infty} + \infty$, and $M_L:= L G_L$.  
We assume that $V_{\rm per} \in C^0_{\rm per}(\Gamma) \cap H^{r-2}_{\rm per}(\Gamma)$ and $W\in C^0(\R^d) \cap H^{r-2}(\R^d)$ for some $r > 2$. Then,
\begin{itemize}
 \item[1.] \bfseries Absence of pollution \normalfont \itshape
\end{itemize}
\begin{equation}\label{eq:nopollutionscint}
\mathop{\lim}_{L\to\infty} \sigma(H_{L,N_L,M_L}) = \sigma(A).
\end{equation}
\begin{itemize}
 \item[2.] \bfseries A priori error estimates
\end{itemize}
Let $\lambda$ be a discrete eigenvalue of $A$ and $\varepsilon >0$ be such that $\sigma(A) \cap (\lambda - \varepsilon, \lambda +\varepsilon) = \{\lambda\}$. 
Let $\cP := \un_{\left\{\lambda \right\}}(A)$ be the $L^2(\R^d)$-orthogonal spectral projection onto the eigenspace of $A$ associated with $\lambda$, and 
$\fP_L := \un_{(\lambda - \varepsilon/2, \lambda +\varepsilon/2)}(H_{L,N_L,M_L})$ the $L^2_{\rm per}(\Gamma_L)$-orthogonal spectral projection of $H_{L,N_L,M_L}$
 associated with the eigenvalues belonging to the interval $(\lambda-\varepsilon/2, \lambda+\varepsilon/2)$. We finally consider a sequence $(\chi_L)_{L \in \N^\ast}$ of cut-off functions such that
$$
0\leq \chi_L\leq 1 \mbox{ on  } \R^d, \; \chi_L  = 1 \mbox{ on } \Gamma_L,
\; \mbox{\rm Supp}(\chi_L)\subset (L+\sqrt{L})\Gamma, \; \|\nabla \chi_L\|_{L^{\infty}} \leq c,
$$
for some constant $c\in\R_+$ independent of $L\in\N^*$. 

Then, $\mbox{\rm Ran}(\cP)\subset H^r(\R^d)$, and there exists $C,\delta>0$ such that for $L$ large enough,
\begin{equation}\label{eq:dimequalityscint}
\mbox{\rm Tr}(\cP) = \mbox{\rm Tr}(\fP_L),
\end{equation}
\begin{equation}\label{eq:proj1scint}
 \mathop{\sup}_{\psi\in{\rm Ran}(\cP), \; \|\psi\|_{L^2(\R^d)} = 1} \mathop{\inf}_{u_L \in {\rm Ran}(\fP_L)} \|\psi - \chi_Lu_L\|_{H^1(\R^d)} \leq  C \left(\epsilon_{1}(L) + \epsilon_2(L) \right),
\end{equation}  
\begin{equation}\label{eq:proj2scint}
 \mathop{\sup}_{u_L\in {\rm Ran}(\fP_L), \; \|u_L\|_{L^2_{\rm per}(\Gamma_L)} = 1} \mathop{\inf}_{\psi \in {\rm Ran}(\cP)} \|\psi - \chi_Lu_L\|_{H^1(\R^d)} \leq  C \left(\epsilon_{1}(L) + \epsilon_2(L) \right),
\end{equation}
\begin{equation}\label{eq:eigenvaluescint}
\mathop{\max}_{\lambda_L \in \sigma(H_{L,N_L})\cap (\lambda-\varepsilon/2, \lambda+\varepsilon/2)} |\lambda_L - \lambda| \leq  C \left(\epsilon_{1}(L)^2 + \epsilon_2(L) \right),
\end{equation}  
where 
$$
\epsilon_1(L):=  e^{-\delta L} + \left(\frac{L}{N_L}\right)^{r-1} \quad \mbox{and} \quad \epsilon_2(L):=  \left(\frac{L}{M_L}\right)^{r-2} + \|W\|_{L^{\infty}\left(\R^d \setminus \Gamma_{L-1}\right)} \left(e^{-\delta L} + \left(\frac{L}{N_L}\right)^{r}\right).
$$
\end{theorem}

\subsection{Formulation in terms of non-consistent approximations}\label{sec:supercellapprox}

The supercell method can be rewritten as a non-consistent approximation of the operator $A$ (in the sense introduced in Section~\ref{sec:approximation}), based on the approximation spaces $(X_L)_{L\in\N^*}$ and the symmetric bilinear forms 
$(a_L)_{L\in\N^*}$, $(\at_L)_{L\in\N^*}$, and $(m_L)_{L\in\N^*}$ defined for all $L\in\N^*$ by
$$
 X_L  :=  \{ \chi_L u_L, \; u_L \in Y_{L,N_L}\} \subset H^1(\R^d),
$$
and
\begin{eqnarray*}
\forall \phi,\psi\in H^1(\R^d), & \quad &  a_L(\phi,\psi) :=   \int_{\Gamma_L} \nabla \phi \cdot \nabla \psi + \int_{\Gamma_L} \cI_{L,M_L}(V_{\rm per} + \widetilde{W}_L)\phi\psi,\\
&&   \at_L(\phi,\psi) := \int_{\Gamma_L} \nabla \phi \cdot \nabla \psi + \int_{\Gamma_L} (V_{\rm per} + W)\phi\psi,\\
&& m_L(\phi,\psi) := \int_{\Gamma_L} \phi\psi,
\end{eqnarray*}
where we recall that $(\chi_L)_{L \in \N^\ast}$ is a sequence of cut-off functions satisfying (\ref{eq:prop_chiL}).
It is easily checked that for all $L\in\N^*$, $m_L(\cdot, \cdot)$ defines a scalar product on $X_L$.

\medskip

Let us introduce, for each $L\in\N^*$, the unitary operator 
$$
\begin{array}{ccc}
 j_L \; : \;  \left( Y_{L,N_L}, \langle \cdot, \cdot \rangle_{L^2_{\rm per}(\Gamma_L)}\right) & \rightarrow & \left( X_L, m_L(\cdot, \cdot)\right),\\
u_L & \mapsto & \chi_Lu_L.\\
\end{array}
$$
Its adjoint (and inverse) $j_L^*$ is given by: $\forall \phi_L \in X_L$, $j_L^*(\phi_L) = u_L$ where $u_L$ is the $L\cR$-periodic extension of $\phi_L|_{\Gamma_L}$. The supercell problems
$$
\left\{
\begin{array}{l}
 \mbox{find }(\lambda_L,u_L)\in \R\times Y_{L,N_L} \mbox{ such that }\|u_L\|_{L^2_{\rm per}(\Gamma_L)} = 1,\\
\forall v_L \in Y_{L,N_L}, \quad \ah_L(u_L, v_L) = \lambda_L \mh_L( u_L, v_L),
\end{array}
\right.
$$
and
$$
\left\{
\begin{array}{l}
 \mbox{find }(\lambda_L,u_L)\in \R\times  Y_{L,N_L} \mbox{ such that }\|u_L\|_{L^2_{\rm per}(\Gamma_L)} = 1,\\
\forall v_L \in Y_{L,N_L}, \quad \ah_{L,M_L}(u_L, v_L) = \lambda_L \mh_L( u_L, v_L),
\end{array}
\right.
$$
are then respectively equivalent, through the change of variable $\psi_L=j_Lu_L$, to the generalized eigenproblems
$$
\left\{
\begin{array}{l}
 \mbox{find }(\lambda_L,\psi_L)\in \R\times X_L \mbox{ such that } m_L(\psi_L,\psi_L) = 1 \mbox{ and }\\
\forall \phi_L  \in X_L, \quad \at_L(\psi_L, \phi_L) = \lambda_L m_L(\psi_L, \phi_L),
\end{array}
\right.
$$
and
$$
\left\{
\begin{array}{l}
 \mbox{find }(\lambda_L,\psi_L)\in \R\times X_L \mbox{ such that } m_L(\psi_L,\psi_L) = 1 \mbox{ and }\\
\forall \phi_L  \in X_L, \quad a_L(\psi_L, \phi_L) = \lambda_L m_L(\psi_L, \phi_L).
\end{array}
\right.
$$

Thus, considering the supercell method with exact and numerical integrations is equivalent to considering the non-consistent but conforming approximations $(\cT_L)_{L\in\N^*}$ and $(\widetilde{\cT}_L)_{L\in\N^*}$ 
respectively defined by
$$
\widetilde{\cT}_L = (X_L, \widetilde{a}_L, m_L) \quad \mbox{and} \quad  \cT_L = (X_L, a_L, m_L).
$$
Taking the same notation as in Section~\ref{sec:theory}, it holds that $\widetilde{A}_L = j_L H_{L,N_L} j_L^*$ and $A_L = j_L H_{L,N_L,M_L} j_L^*$ so that
$\sigma(\widetilde{A}_L) = \sigma(H_{L,N_L})$, $\sigma(A_L) = \sigma(H_{L,N_L,M_L} )$ and, in both cases, $\cP_L = i_{X_L} j_L \mathfrak{P}_L j_L^* i_{X_L}^*$. 
The following section is devoted to the proof of Theorems~\ref{th:supercell} and~\ref{th:supercellint}, which are in fact corollaries of Theorem~\ref{th:thmain}. 
We will first check that all the assumptions of Theorem~\ref{th:thmain} are satisfied for the approximations $(\widetilde{\cT}_L)_{L\in\N^*}$ 
and $({\cT}_L)_{L\in\N^*}$, and then derive more explicit expressions of the right hand sides of (\ref{eq:proj1}), (\ref{eq:proj2}) and (\ref{eq:eigenvalue}) in 
terms of $L$, $N_L$ and $M_L$. 
  
We prove in Section~\ref{sec:assA14} that the supercell method with exact integration satisfies assumptions (A1)-(A4). In Section~\ref{sec:nopollproof}, we prove (\ref{eq:nopollutionsc}) and (\ref{eq:dimequalitysc}), which imply that this method also satisfies assumptions (B1) and (B2) for any discrete eigenvalue $\lambda$ of the operator $A$. 
Estimating the terms involved in estimates (\ref{eq:proj1}), (\ref{eq:proj2}) and (\ref{eq:eigenvalue}) will then lead to estimates (\ref{eq:proj1sc}), (\ref{eq:proj2sc}) and (\ref{eq:eigenvaluesc}) and conclude the proof of Theorem~\ref{th:supercell}. Section~\ref{sec:proofint} is devoted to the proof of Theorem~\ref{th:supercellint}, in which numerical integration errors are taken into account.

\section{Proof of Theorem~\ref{th:supercell} and Theorem~\ref{th:supercellint}}\label{sec:supercellproof}

In the sequel, $C$ will denote an arbitrary constant independent on $L\in\N^*$ which may vary along the calculations.

\subsection{Proof of (A1)-(A4) for $\widetilde{\cT}_L = (X_L, \widetilde{a}_L, m_L)$}\label{sec:assA14}

\noindent
{\bf Proof of (A1):} Let us prove that 
$$
\forall \phi \in H^1(\R^d), \; \mathop{\inf}_{\phi_L\in X_L} \left\|  \phi - \phi_L\right\|_{H^1(\R^d)} \mathop{\longrightarrow}_{L\to\infty} 0.
$$
Let $\phi\in H^1(\R^d)$ and $\varepsilon>0$. Since $C^{\infty}_c(\R^d)$ is dense in $H^1(\R^d)$, there exists $\eta\in C^{\infty}_c(\R^d)$ such that 
$\|\phi - \eta\|_{H^1(\R^d)}\leq \varepsilon$. 
Let $L_0\in\N^*$ be such that $\mbox{\rm Supp}(\eta) \subset (L_0 - \sqrt{L_0})\Gamma$. For all $L\geq L_0$, if $\eta_L$ denotes the $L{\cal R}$-periodic extension of $\eta|_{\Gamma_L}$, we infer from (\ref{eq:Jackson}) that
$$
\|\eta_L - \Pi_{L,N_L}\eta_L\|_{H^1_{\rm per}(\Gamma_L)} \leq C\frac{L}{N_L}\|\eta_L\|_{H^2_{\rm per}(\Gamma_L)} = C\frac{L}{N_L}\|\eta\|_{H^2(\R^d)}\mathop{\longrightarrow}_{L\to\infty} 0,
$$
with $C\in\R_+$ independent of $L$. Let us then consider the sequence $(\phi_L)_{L\in\N^*}$ defined as $\phi_L:= \chi_L \Pi_{L,N_L}\eta_L\in X_L$ for all $L\in\N^*$, for which
\begin{eqnarray*}
 \|\phi - \phi_L\|_{H^1(\R^d)} & \leq & \|\phi - \eta\|_{H^1(\R^d)} + \|\eta - \chi_L \Pi_{L,N_L}\eta_L\|_{H^1(\R^d)},\\
& \leq & \varepsilon + \|\eta_L - \Pi_{L,N_L}\eta_L\|_{H^1_{\rm per}(\Gamma_L)} + \|\chi_L \Pi_{L,N_L}\eta_L\|_{H^1((L+\sqrt{L})\Gamma \setminus\Gamma_L)}.\\
\end{eqnarray*}
Furthermore, since $0 \le \chi_L \le 1$, and $\eta_L = 0$ on $(L+\sqrt{L})\Gamma \setminus \Gamma_L$, it holds
\begin{eqnarray*}
\|\chi_L\Pi_{L,N_L}\eta_L\|^2_{H^1((L+\sqrt{L})\Gamma \setminus \Gamma_L)} &\leq& \left\| \chi_L\Pi_{L,N_L}\eta_L \right\|^2_{L^2((L+\sqrt{L})\Gamma \setminus \Gamma_L) } + 2\left\| \nabla \chi_L \Pi_{L,N_L} \eta_L\right\|^2_{L^2((L+\sqrt{L})\Gamma \setminus \Gamma_L)} \\ && + 2\left\|\chi_L \nabla (\Pi_{L,N_L}\eta_L) \right\|^2_{L^2((L+\sqrt{L})\Gamma \setminus \Gamma_L)}\\
&\leq& \left\| \Pi_{L,N_L}\eta_L - \eta_L \right\|^2_{L^2((L+\sqrt{L})\Gamma \setminus \Gamma_L) } + \left\| \nabla (\Pi_{L,N_L}\eta_L) - \nabla \eta_L \right\|^2_{L^2((L+\sqrt{L})\Gamma \setminus \Gamma_L)}  \\ && + \|\nabla \chi_L\|_{L^{\infty}(\R^d)} \left\| \Pi_{L,N_L} \eta_L - \eta_L\right\|^2_{L^2((L+\sqrt{L})\Gamma \setminus \Gamma_L)} \\
&\leq& 3^d(4 + \|\nabla \chi_L\|_{L^{\infty}(\R^d)}) \left\| \Pi_{L,N_L}\eta_L - \eta_L \right\|^2_{ H^1_{\rm per}(\Gamma_L) } \quad \mathop{\longrightarrow}_{L\to \infty} 0.
\end{eqnarray*}
Hence the result.

\medskip

\noindent
{\bf Proof of (A2):} Let $\phi_L, \psi_L \in X_L$, and $u_L, v_L\in Y_{L,N_L}\subset H^1_{\rm per}(\Gamma_L)$ such that $\phi_L = \chi_L u_L$ and $\psi_L = \chi_L v_L$. It holds
$$
 \int_{\R^d} |\phi_L|^2  =  \int_{\Gamma_{3L}} |\phi_L|^2 =  \int_{\Gamma_{3L}} \chi_L^2 |u_L|^2 \leq  3^d \int_{\Gamma_L}|u_L|^2 =  3^d\int_{\Gamma_L} |\phi_L|^2 \le 3^d \int_{\R^d} |\phi_L|^2.
$$
Therefore, 
\begin{equation} \label{eq:equivL2}
 \frac{1}{3^d}\|\phi_L\|_{L^2(\R^d)}^2 \leq m_L(\phi_L, \phi_L) \leq \|\phi_L\|_{L^2(\R^d)}.
\end{equation}
Besides, 
\begin{eqnarray*}
|\widetilde{a}_L(\phi_L,\psi_L)| &=& \left| \int_{\Gamma_L}\nabla \phi_L\cdot \nabla \psi_L + \int_{\Gamma_L}(V_{\rm per} +W)\phi_L \psi_L\right|\\
&\leq&  (1 + \|W\|_{L^{\infty}(\R^d)}) \|\phi_L\|_{H^1(\R^d)}\|\psi_L\|_{H^1(\R^d)} + \|V_{\rm per}u_L\|_{L^2_{\rm per}(\Gamma_L)}\|v_L\|_{L^2_{\rm per}(\Gamma_L)}\\
&\leq&  ( 1 + \|W\|_{L^{\infty}(\R^d)} + \|V_{\rm per}\|_{\cZ_{\rm per}(\Gamma)}) \|\phi_L\|_{H^1(\R^d)} \|\psi_L\|_{H^1(\R^d)}.
\end{eqnarray*} 
Thus, assumption (A2) is satisfied.

\medskip

\noindent
{\bf Proof of (A3):} For all $\alpha>0$ arbitrarily small, there exists a constant $C_{\alpha}$ such that for all $\phi\in H^1(\R^d)$, 
\begin{equation}\label{eq:semibounded}
\int_{\R^d} |V_{\rm per}||\phi|^2 \leq \alpha \int_{\R^d} |\nabla \phi|^2 + C_{\alpha} \int_{\R^d} |\phi|^2.
\end{equation}
Besides, for all $\phi_L \in X_L$, if $\phi_L = \chi_L u_L$ with $u_L \in Y_{L,N_L}$, it holds that
\begin{eqnarray*}
 \int_{\R^d} |\nabla \phi_L|^2 &\leq & 2 \int_{(L+\sqrt{L})\Gamma} |\nabla \chi_L u_L|^2 + |\chi_L \nabla u_L|^2 \\
& \leq & 2 \times  3^d \left( \|\nabla \chi_L\|_{L^{\infty}(\R^d)} \int_{\Gamma_L} |u_L|^2 + \int_{\Gamma_L} |\nabla u_L|^2\right),\\
\end{eqnarray*}
which, together with (\ref{eq:equivL2}), yields that, for $L$ large enough 
\begin{equation}\label{eq:eqgrad}
 \int_{\R^d} |\phi_L|^2 + \int_{\R^d} |\nabla \phi_L|^2 \leq 3^{d+1} \left( \int_{\Gamma_L} |\phi_L|^2 + \int_{\Gamma_L} |\nabla \phi_L|^2\right).
\end{equation}
Using (\ref{eq:semibounded}) and (\ref{eq:eqgrad}), we obtain that for all $\alpha >0$ arbitrarily small, there exists $D_\alpha\in \R_+$ such that for all 
$L\in\N^*$ and all $\phi_L \in X_L$, 
\begin{eqnarray*}
\int_{\Gamma_L} (V_{\rm per} + W)|\phi_L|^2 & \leq & \int_{\R^d} \left(|V_{\rm per}| + |W|\right)|\phi_L|^2\\
& \leq & \alpha \int_{\Gamma_L} |\nabla \phi_L|^2 + D_\alpha \int_{\Gamma_L} |\phi_L|^2.\\
\end{eqnarray*}
This last inequality implies that there exists $\beta>0$ independent on $L\in\N^*$ such that for all $\phi_L\in X_L$, 
$$
\|\phi_L\|_{H^1(\R^d)}^2 \leq 3^{d+1} \|\phi_L\|_{H^1(\Gamma_L)}^2 \leq \beta\left(|\widetilde{a}_L(\phi_L, \phi_L)| +  m_L(\phi_L, \phi_L)\right).
$$
Thus, for all $\mu\in \C$, it holds that
\begin{eqnarray*}
&&\mathop{\inf}_{\psi_L\in X_L} \mathop{\sup}_{\phi_L\in X_L} \frac{|(\at_L-\mu m_L)(\phi_L, \psi_L)|}{\|\phi_L\|_{H^1(\R^d)} \|\psi_L\|_{H^1(\R^d)}} \\
&&\geq 
\frac{1}{\beta}\mathop{\inf}_{\psi_L\in X_L} \mathop{\sup}_{\phi_L\in X_L} \frac{|(\at_L-\mu m_L)(\phi_L, \psi_L)|}{\left( |\widetilde{a}_L(\phi_L, \phi_L)| + m_L(\phi_L, \phi_L) \right)^{1/2} 
\left( |\widetilde{a}_L(\psi_L, \psi_L)| + m_L(\psi_L, \psi_L) \right)^{1/2}}.
\end{eqnarray*}
Let $(\zeta_L^i)_{1\leq i \leq \mbox{\rm dim}(X_L)}$ be an $m_L$-orthonormal basis of $X_L$, such that for all $1\leq i\leq \mbox{\rm dim}(X_L)$, 
$$
H_{L,N_L}j_L^* \zeta_L^{(i)} = \nu_L^i j_L^*\zeta_L^i, \; 1\leq i \leq \mbox{\rm dim}(X_L),
$$
where $\{ \nu_L^i, \; 1\leq i \leq \mbox{\rm dim}(X_L)\} = \sigma(H_{L,N_L})$. Then, any $\phi_L\in X_L$ can be expanded in the basis $(\zeta_L^i)_{1\leq i \leq \mbox{\rm dim}(X_L)}$:
$$
\phi_L = \sum_{i=1}^{\mbox{\rm dim}(X_L)} c_i \zeta_L^i, \; c_i \in \R, \; 1\leq i \leq \mbox{\rm dim}(X_L),
$$
and it holds that $|\widetilde{a}_L(\phi_L, \phi_L)| + m_L(\phi_L, \phi_L) \leq \sum_{i=1}^{\mbox{\rm dim}(X_L)} |c_i|^2 (1+|\nu_L^i|)$. Considering 
$$
\psi_L:= \sum_{i=1}^{\mbox{\rm dim}(X_L)} \mbox{\rm sgn}(\nu_L^i -\mu)c_i \zeta_L^i, 
$$
we obtain that
\begin{equation}\label{eq:eqcentr}
\mathop{\inf}_{\psi_L\in X_L} \mathop{\sup}_{\phi_L\in X_L} \frac{|(\at_L-\mu m_L)(\phi_L, \psi_L)|}{\|\phi_L\|_{H^1(\R^d)} \|\psi_L\|_{H^1(\R^d)}} \geq \frac{1}{\beta} \mathop{\inf}_{\nu_L \in \sigma(H_{L,N_L})} \frac{|\nu_L - \mu|}{1+|\nu_L|}.
\end{equation}
Since (\ref{eq:eqcentr}) holds for any $\mu\in \C$, this implies that for any compact subset $K\subset \C$, there exists a constant $c_K>0$ such that for all $L\in\N^*$ and all $\mu\in K$, 
$$
\mathop{\inf}_{\psi_L\in X_L} \mathop{\sup}_{\phi_L\in X_L} \frac{|(\at_L-\mu m_L)(\phi_L, \psi_L)|}{\|\phi_L\|_{H^1(\R^d)} \|\psi_L\|_{H^1(\R^d)}} \geq c_K\min \left( 1, \mbox{\rm dist}(\mu, \sigma(H_{L,N_L}))\right).
$$
Thus, condition (A3'), and condition (A3), hold for the approximation $(\widetilde{\cT}_L)_{L\in \N^*}$. 

\medskip

\noindent
{\bf Proof of (A4):} For all $\phi\in H^1(\R^d)$, we denote by
$$
r_L^m(\phi): = \left( \int_{\R^d\setminus \Gamma_{L-1}} |\phi|^2\right)^{1/2} \leq \|\phi\|_{L^2(\R^d)},
$$
and
$$
r_L^a(\phi) := \left( \int_{\R^d\setminus \Gamma_{L-1}}|\phi|^2+ |\nabla \phi|^2 \right)^{1/2}\leq \|\phi\|_{H^1(\R^d)}.
$$
Then, $r_L^m$ and $r_L^a$ are seminorms on $H^1(\R^d)$ such that for all $\phi\in H^1(\R^d)$, $\dps r_L^m(\phi)\mathop{\longrightarrow}_{L\to\infty} 0$ and 
$\dps r_L^a(\phi) \mathop{\longrightarrow}_{L\to\infty} 0$. 
For all $\phi,\psi \in H^1(\R^d)$, it holds 
$$
 |(m-m_L)(\phi, \psi)|  =  \left| \int_{\R^d \setminus \Gamma_L} \phi \psi \right| \leq  r_L^m(\phi)r_L^m(\psi).  
$$
Let $(\omega_L)_{L\in\N^*}$ be a sequence of $C^{\infty}$ cut-off functions such that for all $L\in\N^*$, $0\leq \omega_L \leq 1$, $\omega_L = 1$ on $\R^d\setminus \Gamma_L$, 
$\omega_L = 0$ on $\Gamma_{L-1}$ and the sequence $(\|\nabla \omega_L\|_{L^\infty(\R^d)})_{L\in\N^*}$ is uniformly bounded in $L\in\N^*$. 
Then, for all $\phi,\psi\in H^1(\R^d)$, 
\begin{eqnarray*}
 |(a-\at_L)(\phi,\psi)| & = & \left| \int_{\R^d\setminus \Gamma_L} \nabla \phi \cdot \nabla \psi + \int_{\R^d\setminus \Gamma_L} (V_{\rm per} + W) \phi\psi\right|\\
& \leq & (1 + \|W\|_{L^{\infty}}) r_L^a(\phi) r_L^a(\psi) + \int_{\R^d} |V_{\rm per} \omega_L\phi \omega_L\psi |\\
& \leq & (1 + \|W\|_{L^{\infty}}) r_L^a(\phi) r_L^a(\psi) + \left(\int_{\R^d} |V_{\rm per}| |\omega_L\phi|^2\right)^{1/2} \left(\int_{\R^d} |V_{\rm per}| |\omega_L\psi|^2\right)^{1/2}.\\
\end{eqnarray*}
Using (\ref{eq:semibounded}),
$$
\int_{\R^d} |V_{\rm per}| |\omega_L\phi|^2 \leq  \frac{1}{2} \int_{\R^d} |\nabla (\omega_L \phi)|^2 + C \int_{\R^d} |\omega_L \phi|^2 \leq C r_L^a(\phi)^2.
$$
Thus, there exists $\kappa \in \R_+$ independent on $L\in\N^*$ such that
$$
|(a-\at_L)(\phi,\psi)|  \leq  \kappa r_L^a(\phi) r_L^a(\psi).
$$

\subsection{Absence of pollution} \label{sec:nopollproof}

\begin{proposition}\label{prop:nopollution}
It holds
\begin{equation}\label{eq:limeigset}
\sigma(A) = \mathop{\lim}_{L\to\infty} \sigma(H_{L,N_L}).
\end{equation}
Besides, for any discrete eigenvalue $\lambda$ of the operator $A$ and for all $\varepsilon>0$ such that $(\lambda -\varepsilon, \lambda + \varepsilon) \cap \sigma(A) = \{\lambda\}$, we have, for $L$ large enough,
\begin{equation}\label{eq:dimequality}
\mbox{\rm Tr}(\fP_L) = \mbox{\rm Tr}(\cP),
\end{equation}
where $\cP := \un_{\left\{\lambda\right\}}(A)$ and $\fP_L := \un_{(\lambda-\varepsilon/2, \lambda+\varepsilon/2)}(H_{L,N_L})$. 
\end{proposition}

Let us notice that (\ref{eq:limeigset}) implies that (B1) is satisfied for any discrete eigenvalue of $A$, and that (\ref{eq:dimequality}) is nothing but a reformulation of (B2). We refer to \cite[Theorem~3.1]{CEM} for a proof of (\ref{eq:limeigset}). 

\begin{proof}[Proof of (\ref{eq:dimequality})]
If follows from  (\ref{eq:limeigset}) that (B1) is satisfied and therefore that for $n$ large enough, $\mbox{\rm Tr}(\fP_{L_n}) \ge \mbox{\rm Tr}(\cP)$. Let us assume that there exists an increasing sequence $(L_k)_{k\in\N^*}$ of integers such that 
$$
\mbox{\rm Tr}(\fP_{L_k}) > q:=\mbox{\rm Tr}(\cP).
$$
For all $k\in\N$, let $(\zeta_{L_k}^{(i)})_{1\leq i \leq q+1}$ be an $L^2_{\rm per}(\Gamma_{L_k})$-orthonormal family of vectors of $Y_{L_k, N_{L_k}}$ such that for all $1\leq i \leq q+1$, 
$$
H_{L_k, N_{L_k}}\zeta_{L_k}^{(i)} = \lambda_{L_k}^{(i)}\zeta_{L_k}^{(i)} \quad \mbox{with }\lambda_{L_k}^{(i)}\in (\lambda - \varepsilon/2, \lambda +\varepsilon/2).
$$
Then, for all $k\in\N$, $(\chi_{L_k}\zeta_{L_k}^{(i)})_{1\leq i \leq q+1}$ forms a free family of $X_{L_k}$ and there exists 
$g_k\in \mbox{\rm Span}(\zeta_{L_k}^{(i)})_{1\leq i \leq q+1}$ such that $\|g_k\|_{L^2_{\rm per}(\Gamma_{L_k})} = 1$ and
$$
\widetilde{g}_k := \chi_{L_k}g_k \in \mbox{\rm Ker}(\cP).
$$
Reasoning as above, it can be easily checked that $\left( \|g_k\|_{H^1_{\rm per}(\Gamma_{L_k})}\right)_{k\in\N^*}$ is bounded, which implies that 
$ \left(\| \widetilde{g}_k\|_{H^1(\R^d)}\right)_{k\in\N^*}$ is bounded as well. Thus, up to the extraction of a subsequence, there exists $g\in H^1(\R^d) \cap \mbox{\rm Ker}(\cP)$ 
such that 
$\dps \widetilde{g}_k \mathop{\wlim}_{k\to\infty} g$ in $H^1(\R^d)$ and $\dps \widetilde{g}_k \mathop{\longrightarrow}_{k\to\infty} g$ in $L^2_{\rm loc}(\R^d)$. Since 
$\widetilde{g}_k=\chi_{L_k}g_k$, this also implies that
$$
g_k \mathop{\longrightarrow}_{k\to\infty} g \quad \mbox{strongly in }L^2_{\rm loc}(\R^d),
$$
which readily leads to 
$$
\left(H_{L_k,N_{L_k}} -\lambda\right) g_k \mathop{\longrightarrow}_{k\to\infty} -\Delta g + (V_{\rm per} + W -\lambda)g \quad \mbox{ in }\cD'(\R^d).
$$
Besides, since $g_k \in \mbox{\rm Ran}(\fP_{L_k})$ and $\dps \mathop{\lim}_{k\to\infty} \sigma(H_{L_k,N_{L_k}}) = \sigma(H)$, we have,
$$
\left\|  \left(H_{L_k,N_{L_k}} -\lambda\right) g_k \right\|_{L^2_{\rm per}(\Gamma_{L_k})} \mathop{\longrightarrow}_{k\to\infty} 0,
$$
which, in turn, implies that
$$
\left(H_{L_k,N_{L_k}} -\lambda\right) g_k \mathop{\longrightarrow}_{k\to\infty} 0 \quad \mbox{ in }\cD'(\R^d).
$$
Therefore,
$$
 -\Delta g + (V_{\rm per} + W -\lambda)g =0.
$$
Consequently, $g\in \mbox{\rm Ker}(\cP) \cap \mbox{\rm Ran}(\cP) = \{0\}$. Using similar arguments as those used in the proof of \cite[Theorem~3.1]{CEM}, we infer from the fact that $(g_k)_{k \in \N}$ strongly converges to $0$ in $L^2_{\rm loc}(\R^d)$ that $\left( \frac{\widetilde{g}_k}{\|\widetilde{g}_k\|_{L^2(\R^d)}}\right)_{k\in\N}$ is 
a Weyl sequence for $A^0=-\Delta + V_{\rm per}$ associated with $\lambda$, which contradicts the fact that $\lambda \notin \sigma(A^0)$.
\end{proof}

\subsection{Proof of Theorem~\ref{th:supercell}}

We have proved that the supercell method with planewave discretization and exact integration satisfies assumptions (A1)-(A4), and that for each discrete eigenvalue located in a spectral gap of $A$, assumptions (B1) and (B2) are satisfied. Thus, Theorem~\ref{th:thmain} can be applied and there exists $C\in \R_+$ such that for $L$ large enough, 
\begin{eqnarray*}
  \mbox{\rm Tr}(\cP_L) & = & \mbox{\rm Tr}(\cP) = \mbox{\rm Tr}(\fP_L),\\
 \|(\cP-\cP_L)\cP\|_{\cL(L^2(\R^d),H^1(\R^d))} & \leq & C  \left(\left\| \left( 1 - \Pi_{X_L}^{H^1(\R^d)}\right) \cP \right\|_{\cL(L^2(\R^d),H^1(\R^d))} + \cR_L^a + \cR_L^m\right),\\
 \|(\cP-\cP_L)\cP_L\|_{\cL(L^2(\R^d),H^1(\R^d))}&  \leq & C\left( \left\| \left(1-\Pi_{X_L}^{H^1(\R^d)}\right) \cP \right\|_{\cL(L^2(\R^d),H^1(\R^d))} + \cR_L^a + \cR_L^m\right),\\
\mathop{\max}_{\lambda_L \in \sigma(H_{L,N_L})\cap (\lambda - \varepsilon/2, \lambda+\varepsilon/2)} |\lambda_L - \lambda| & \leq &C\left( \left\| \left( 1 - \Pi_{X_L}^{H^1(\R^d)} \right) \cP \right\|_{\cL(L^2(\R^d),H^1(\R^d))} + \cR_L^a + \cR_L^m \right)^2,
\end{eqnarray*}
where $\cP_L:= i_{X_L} j_L \fP_L j_L^* i_{X_L}^*$ and 
\begin{eqnarray*}
 \cR_L^m &:= & \mathop{\sup}_{\psi\in {\rm Ran}(\cP), \; \|\psi\|_{L^2(\R^d)} = 1} r_L^m(\psi),\\
 \cR_L^a &:= & \mathop{\sup}_{\psi\in {\rm Ran}(\cP), \; \|\psi\|_{L^2(\R^d)} = 1} r_L^a(\psi).
\end{eqnarray*}
Since we have
$$
 \mathop{\sup}_{\psi\in {\rm Ran}(\cP), \; \|\psi\|_{L^2(\R^d)} = 1} \mathop{\inf}_{u_L \in {\rm Ran}(\fP_L)} \|\psi - \chi_Lu_L\|_{H^1(\R^d)}  \leq  \|(\cP-\cP_L)\cP\|_{\cL(L^2(\R^d),H^1(\R^d))},
$$
and 
$$
 \mathop{\sup}_{u_L\in {\rm Ran}(\fP_L), \; \|u_L\|_{L^2_{\rm per}(\Gamma_L)} = 1} \mathop{\inf}_{\psi \in {\rm Ran}(\cP)} \|\psi - \chi_Lu_L\|_{H^1(\R^d)} \leq \|(\cP-\cP_L)\cP_L\|_{\cL(L^2(\R^d),H^1(\R^d))},
$$
it just remains to prove that there exists $\delta>0$ independent on $L$ such that
$$
\left\| \left( 1 - \Pi_{X_L}^{H^1(\R^d)}\right) \cP \right\|_{\cL(L^2(\R^d),H^1(\R^d))} + \cR_L^a + \cR_L^m \leq C\left( e^{-\delta L} + \left(\frac{N_L}{L}\right)^{r-1} \right).
$$
This estimate is based on exponential decay results for the bound states of Schr\"odinger operators~\cite{Simon}. A real-valued function $V$ on $\R^d$ is said to lie in the class $K_d$ if and only if
$$
\begin{array}{ll}
\mbox{ if }d\geq 3, & \dps \quad \lim_{\alpha \downarrow 0}  \sup_{x\in\R^d} \int_{|x-y|\leq \alpha} \frac{|V(y)|}{|x-y|^{d-2}}\,dy =0;\\
\mbox{ if } d=2, & \dps \quad \lim_{\alpha \downarrow 0}  \sup_{x\in\R^d} \int_{|x-y|\leq \alpha} |V(y)| \ln\left(|x-y|^{-1}\right)\,dy =0;\\
\mbox{ if } d=1, & \dps \quad \sup_{x\in\R^d} \int_{|x-y|\leq 1} |V(y)|\,dy < \infty.\\
\end{array}
$$
Under our assumptions on $V_{\rm per}$ and $W$, $V = V_{\rm per} + W \in K_d$. It then follows from 
Theorem C.3.4 and Corollary C.2.3 in \cite{Simon} that there exists $C,\delta>0$ such that for all 
$L^2(\R^d)$-normalized $\psi\in \mbox{\rm Ran}(\cP)$, 
\begin{equation}\label{eq:expdecay}
\forall x\in\R^d,\quad |\psi(x)| \leq Ce^{-3\delta|x|} \quad \mbox{and} \quad 
 e^{3\delta |\cdot|} \nabla \psi \in \left(L^2(\R^d)\right)^d.
\end{equation}

For all $L\geq 6$, let $\eta_L\in C^{\infty}_c(\R^d)$ 
such that $0\leq \eta_L \leq 1$, $\eta_L = 1 $ on $\Gamma_{L/2-2}$, $\mbox{\rm Supp}(\eta_L) \subset \Gamma_{L/2-1}$ and all 
its derivative up to the $[r+1]^{\rm st}$ order are bounded in $L^{\infty}(\R^d)$, uniformly in $L\in\N^*$.
Let $\psi\in\mbox{\rm Ran}(\cP)$ such that $\|\psi\|_{L^2(\R^d)} = 1$, $\zeta_L = \eta_L\psi$, and $\widetilde{\zeta}_L$ the $L\cR$-periodic extension of $\zeta_L$. Then, 
$\chi_L \Pi_{L,N_L}\widetilde{\zeta}_L \in X_L$, and it holds
\begin{eqnarray*}
 \|\psi - \chi_L \Pi_{L,N_L}\widetilde{\zeta}_L\|_{H^1(\R^d)} & \leq & \|\psi - \eta_L\psi\|_{H^1(\R^d)} + \|\zeta_L - \chi_L\Pi_{L,N_L}\widetilde{\zeta}_L\|_{H^1(\R^d)}\\
& = & \|\psi - \eta_L\psi\|_{H^1(\R^d)} + \|\chi_L\left(\zeta_L - \Pi_{L,N_L}\widetilde{\zeta}_L\right)\|_{H^1(\R^d)}\\
& \leq & C e^{-\delta L} + C \| \widetilde{\zeta}_L - \Pi_{L,N_L}\widetilde{\zeta}_L\|_{H^1_{\rm per}(\Gamma_L)}\\
& \leq & C e^{-\delta L} + C\left(\frac{N_L}{L}\right)^{r-1} \|\widetilde{\zeta}_L\|_{H^r_{\rm per}(\Gamma_L)}\\
& \leq & C \left( e^{-\delta L} + \left(\frac{N_L}{L}\right)^{r-1} \|\zeta_L\|_{H^r(\R^d)}\right)\\
& \leq & C \left( e^{-\delta L} + \left(\frac{N_L}{L}\right)^{r-1} \|\psi\|_{H^r(\R^d)}\right)\\
& \leq & C \left( e^{-\delta L} + \left(\frac{N_L}{L}\right)^{r-1}\right).
\end{eqnarray*}
This yields the estimate
$$
\left\| \left( 1 - \Pi_{X_L}^{H^1(\R^d)}\right) \cP \right\|_{\cL(L^2(\R^d),H^1(\R^d))}\leq C\left( e^{-\delta L} + \left(\frac{N_L}{L}\right)^{r-1} \right).
$$
The remaining estimate
$$
\cR_L^a + \cR_L^m \leq  C e^{-\delta L},
$$
is a straightforward consequence of (\ref{eq:expdecay}).

\subsection{Proof of Theorem~\ref{th:supercellint}}\label{sec:proofint}

Let us first remark that since $\dps \frac{M_L}{L} = G_L \in \N^\ast$, 
$\cI_{L,M_L}(V_{\rm per}) = \cI_{1,G_L}(V_{\rm per})$ is a $\cR$-periodic function. 
Let $\phi_L, \psi_L \in X_L$ be such that $\phi_L = \chi_L u_L$ and $\psi_L = \chi_L v_L$ with $u_L, v_L \in Y_{L,N_L}$.
Then, we have
$$
\left| \int_{\Gamma_L} (V_{\rm per} - \cI_{L,M_L}(V_{\rm per})) \, \phi_L \psi_L \right| \leq  \left| \int_{\Gamma_L} |V_{\rm per} - \cI_{L,M_L}(V_{\rm per})|  \, u_L^2 \right|^{1/2}\left| \int_{\Gamma_L} |V_{\rm per} - \cI_{L,M_L} (V_{\rm per})| \, v_L^2 \right|^{1/2} .
$$
As 
\begin{eqnarray*} 
\int_{\Gamma_L} |V_{\rm per} - \cI_{L,M_L}(V_{\rm per})|  \, u_L^2 &=&  \sum_{R\in \cR \cap \Gamma_L} \int_{\Gamma} |V_{\rm per} - \cI_{1,G_L}(V_{\rm per})| \,  u_L(\cdot +R)^2 \\
&& \leq  \|V_{\rm per} - \cI_{1,G_L}(V_{\rm per}) \|_{L^2_{\rm per}(\Gamma)} \sum_{R\in \cR \cap \Gamma_L} \|u_L(\cdot +R)\|_{L^4(\Gamma)}^2\\
&& \leq  C \, \|V_{\rm per} - \cI_{1,G_L}(V_{\rm per}) \|_{L^2_{\rm per}(\Gamma)}  \sum_{R\in \cR \cap \Gamma_L} \|u_L(\cdot +R)\|_{H^1(\Gamma)}^2\\
&& =  C \,  \|V_{\rm per} - \cI_{1,G_L}(V_{\rm per}) \|_{L^2_{\rm per}(\Gamma)}  \|u_L\|_{H^1_{\rm per}(\Gamma_L)}^2 \\
&& \le C \,  \|V_{\rm per} - \cI_{1,G_L}(V_{\rm per}) \|_{L^2_{\rm per}(\Gamma)}  \|\phi_L\|_{H^1(\R^d)}^2,
\end{eqnarray*}
we obtain
\begin{equation}\label{eq:ineqVper}
\left| \int_{\Gamma_L} (V_{\rm per} - \cI_{L,M_L}(V_{\rm per})) \phi_L \psi_L\right| \leq C \|V_{\rm per} - \cI_{1,G_L}(V_{\rm per}) \|_{L^2_{\rm per}(\Gamma)} \|\phi_L\|_{H^1(\R^d)} \|\psi_L\|_{H^1(\R^d)},
\end{equation}
for a constant $C$ independent of $L$, with 
\begin{equation}\label{eq:ineqVper2}
\|V_{\rm per} - \cI_{1,G_L}(V_{\rm per}) \|_{L^2_{\rm per}(\Gamma)} \leq C G_L^{-(r-2)} \|V_{ \rm per}\|_{H^{r-2}_{\rm per}(\Gamma)} =  C\left(\frac{L}{M_L}\right)^{r-2}   \|V_{ \rm per}\|_{H^{r-2}_{\rm per}(\Gamma)} \mathop{\longrightarrow}_{L\to\infty} 0.
\end{equation}
Besides, since $W\in C^0(\R^d) \cap H^{r-2}(\R^d)$, 
\begin{eqnarray}
\|\widetilde{W}_L - \cI_{L,M_L}(\widetilde{W}_L) \|_{L^2(\Gamma_L)} & \leq & C\left(\frac{L}{M_L}\right)^{r-2}  \|\widetilde{W}_L \|_{H^{r-2}_{\rm per}(\Gamma_L)}\nonumber \\ 
& \leq & C\left(\frac{L}{M_L}\right)^{r-2}  \|W \|_{H^{r-2}(\R^d)}, \label{eq:ineqW}
\end{eqnarray}
and $\dps \|W - \widetilde{W}_L\|_{L^{\infty}(\R^d)} \leq \|W\|_{L^\infty(\R^d \setminus \Gamma_{L-1})} \mathop{\longrightarrow}_{L\to\infty} 0$. 
Thus,
$$
\mathop{\sup}_{\phi_L\in X_L} \mathop{\sup}_{\psi_L \in X_L} \frac{|(\at_L-a_L)(\phi_L, \psi_L)|}{\|\phi_L\|_{H^1(\R^d)} \|\psi_L\|_{H^1(\R^d)}} \mathop{\longrightarrow}_{L\to \infty} 0.
$$
Together with the results proved in Section~\ref{sec:assA14}, this implies (A2), (A3), (B1) and (B2) are satisfied for ${\cT}_L = (X_L,{a}_L, m_L)$. Assumption (A4) is also satisfied for ${\cT}_L = (X_L,{a}_L, m_L)$, with $\at_L(\cdot, \cdot)$ playing the role of $\at_n(\cdot, \cdot)$ and $\mt_n(\cdot, \cdot)=m_n(\cdot, \cdot)=m_L(\cdot, \cdot)$.
To obtain the estimates (\ref{eq:proj1scint}), (\ref{eq:proj2scint}) and (\ref{eq:eigenvaluescint}), it remains to prove that 
$$
\cS_L^a \leq C \left[ \left( \frac{L}{M_L} \right)^{r-2} + \|W\|_{L^\infty(\R^d \setminus \Gamma_{L-1})} \left( e^{-\delta L} + \left(\frac{L}{N_L}\right)^r \right)\right],
$$
where 
$$
\cS_L^a: = \mathop{\sup}_{\psi \in {\rm Ran}(\cP), \; \|\psi\|_{L^2(\R^d)} = 1} \mathop{\sup}_{\phi_L\in X_L} \frac{\left|(a_L-\at_L)\left(\Pi_{X_L}^{H^1(\R^d)}\psi, \phi_L\right)\right|}{\|\phi_L\|_{H^1(\R^d)}}.
$$
Using (\ref{eq:ineqVper}) and (\ref{eq:ineqVper2}), we already have for all $\psi \in {\rm Ran}(\cP)$ such that $\|\psi\|_{L^2(\R^d)} = 1$,
\begin{equation}\label{eq:Vperterm}
\left| \int_{\Gamma_L} (V_{\rm per} - \cI_{L,M_L}(V_{\rm per}))\left(\Pi_{X_L}^{{H^1(\R^d)}}\psi\right)\phi_L \right| \leq C \left( \frac{L}{M_L}\right)^{r-2} \|\phi_L\|_{H^1(\R^d)}.
\end{equation}
Besides, using (\ref{eq:ineqW}), it holds that
\begin{equation}\label{eq:Wterm1}
\left| \int_{\Gamma_L} (\widetilde{W}_L - \cI_{L,M_L}(\widetilde{W_L}))\left(\Pi_{X_L}^{{H^1(\R^d)}}\psi\right) \phi_L\right| \leq C \left( \frac{L}{M_L}\right)^{r-2} \|\phi_L\|_{H^1(\R^d)}.
\end{equation}
It also follows from (\ref{eq:expdecay}) that
\begin{eqnarray*}
 \left| \int_{\Gamma_L} (\widetilde{W}_L -W) \left(\Pi_{X_L}^{{H^1(\R^d)}}\psi\right)\phi_L \right| & \leq &   
\left| \int_{\Gamma_L\setminus\Gamma_{L-1}} (\widetilde{W}_L -W) \psi \phi_L \right| + \left| \int_{\Gamma_L\setminus\Gamma_{L-1}} (\widetilde{W}_L -W) \left(\psi - \Pi_{X_L}^{H^1(\R^d)}\psi\right) \phi_L \right|\\
& \leq & C\|W\|_{L^\infty(\R^d \setminus \Gamma_{L-1})}\left(  e^{-\delta L} + \left\|\psi - \Pi_{X_L}^{H^1(\R^d)}\psi\right\|_{L^2(\R^d)} \right)\|\phi_L\|_{H^1(\R^d)}.
\end{eqnarray*}
Reasoning as in the proof of (A1) in Section~\ref{sec:assA14}, and using (\ref{eq:expdecay}), we can prove that
$$
  \left\|\psi - \Pi_{X_L}^{H^1(\R^d)}\psi\right\|_{L^2(\R^d)} \leq C \left( e^{-\delta L} + \left(\frac{L}{N_L}\right)^r\right).
$$
Thus, 
\begin{equation}\label{eq:Wterm2}
 \left| \int_{\Gamma_L} (\widetilde{W}_L -W) \left(\Pi_{X_L}^{{H^1(\R^d)}}\psi\right)\phi_L \right| \leq C \left[ \left(\frac{L}{M_L}\right)^{r-2} + \|W\|_{L^\infty(\R^d \setminus \Gamma_{L-1})}\left( e^{-\delta L} + \left(\frac{L}{N_L}\right)^r  \right)\right] \|\phi_L\|_{H^1(\R^d)}.
\end{equation}
Finally, using (\ref{eq:Vperterm}), (\ref{eq:Wterm1}) and (\ref{eq:Wterm2}), we obtain
$$
\cS_L^a \leq C \left[ \left( \frac{L}{M_L} \right)^{r-2} + \|W\|_{L^\infty(\R^d \setminus \Gamma_{L-1})} \left( e^{-\delta L} + \left(\frac{L}{N_L}\right)^r \right)\right],
$$ 
which ends the proof of Theorem~\ref{th:supercellint}.

\section{Numerical results}\label{sec:numerical}

In this section, we present some numerical results obtained with the software Scilab, illustrating the {\it a priori} estimates given in Theorem~\ref{th:supercell} and~\ref{th:supercellint}. These results have been obtained with $d=1$, $V_{\rm per}(x) = |\sin x|$, $W(x) = -2\exp(-|x|)$ and $\Gamma = (-\pi, \pi]$. The particular form of these potentials enables us to compute the mass and stiffness matrices analytically (and therefore with no numerical integration error). 
The operator $A = -\Delta + V_{\rm per} +W$ then possesses a discrete simple eigenvalue $\lambda \approx 1.69$ located in the spectral gap $[\alpha, \beta]$ of the operator $A^0 = -\Delta + V_{\rm per}$ where $\alpha \approx 1.43$ and $\beta \approx 1.84$. The reference values for $\lambda$ and the associated eigenvector (considered in our numerical study as the limits $L,N_L \to \infty$) are obtained with $L_{\rm ref}=40$ and $N_{\rm ref}=1400$. 

\medskip

Figure~1 shows $\sigma(H_{L,N_{\rm ref}})\cap [1,2]$ for $L=6,8,10,12,14,16,18$ and $N_{\rm ref}=1400$. We can see that there is no spectral pollution, as predicted by \cite{CEM} and Proposition~\ref{prop:nopollution}. 

\medskip

\begin{figure}\label{fig:fig1}
\centering
 \includegraphics[width=14cm]{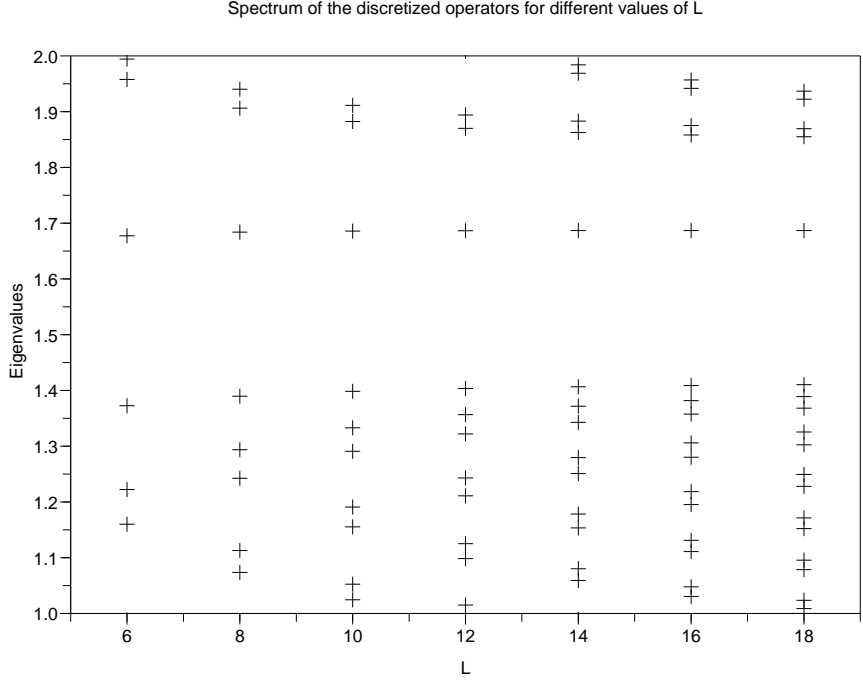}
 % spec_supercell.eps: 0x0 pixel, 300dpi, 0.00x0.00 cm, bb=14 220 581 621
\caption{Spectrum of $H_{L,N_{\rm ref}}$ in the range $[1,2]$ for different values of $L$, with $N_{\rm ref} = 1400$.}
\end{figure}

\medskip

The next series of numerical tests confirms the exponential convergence of the supercell method with respect to the size of the supercell. We have compared 
the eigenvalue closest to $\lambda$ and the associated eigenvector obtained for different values of $L$ ($L=6,8,10,12,14,16,18$) to the reference eigenvalue 
and eigenvector obtained with $L=40$, all these calculations being done with $N_{\rm ref}=1400$. Figure~2 shows the relative errors on the eigenvalue, 
and the square of the $L^2$ and $H^1$ norms of the error on the eigenvector. More precisely, for all $L\in\N^*$, we consider the eigenvector $u_L$ of $H_{L,N_{\rm ref}}$ associated with the eigenvalue $\lambda_L$ of $H_{L,N_{\rm ref}}$ closest to $1.69$, and set $\phi_L = \chi_L u_L$, where $\chi_L$ is the unique $C^2$ function defined by
$\chi_L = 1$ on $[-\pi L, \pi L]$, $\chi_L=0$ on $\R \setminus [-\pi (L+\sqrt L), \pi (L+\sqrt L)]$, and $\chi_L$ is a sixth degree polynomial on $[-\pi (L+\sqrt L),-\pi L]$ and on $[\pi L, \pi (L+\sqrt L)]$. Figure~2 shows the decay rate of $\log_{10}\left(\frac{|\lambda_L -  \lambda_{L_{\rm ref}}|}{\lambda_{L_{\rm ref}}}\right)$, 
$\log_{10}\left(\frac{\|\phi_L -  \phi_{L_{\rm ref}}\|^2_{L^2(\R)}}{\|\phi_{L_{\rm ref}}\|_{L^2(\R)}^2}\right)$ and $\log_{10}\left(\frac{\|\phi_L -  \phi_{L_{\rm ref}}\|^2_{H^1(\R)}}{\|\phi_{L_{\rm ref}}\|_{H^1(\R)}^2}\right)$.
These numerical results show the exponential decay of the error as a function of $L$, as well as the doubling of the convergence rate of the eigenvalue with respect to the convergence rate of the eigenvector.   

\begin{figure}\label{fig:fig2}
\centering
 \includegraphics[width=14cm]{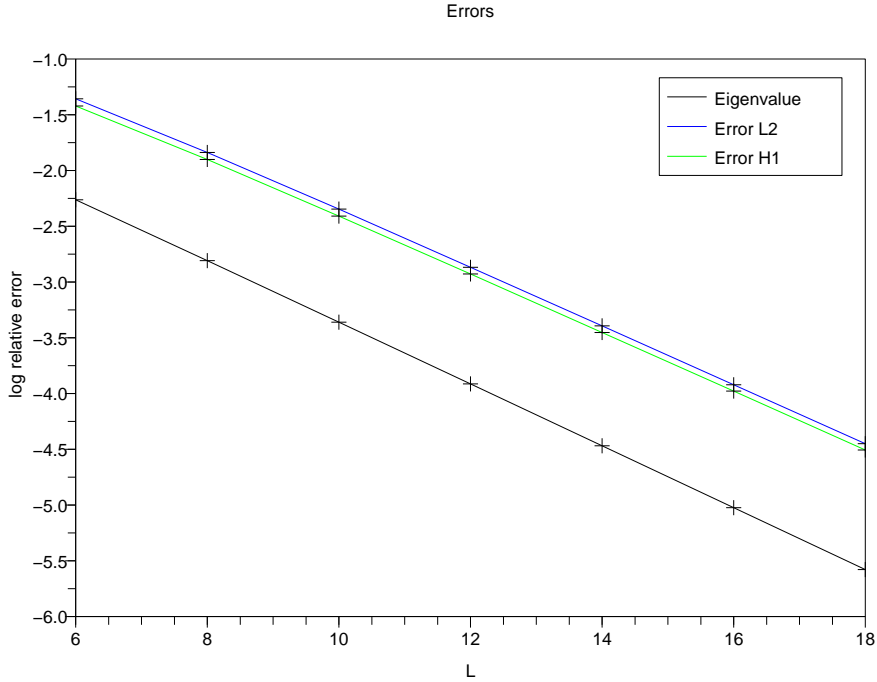}
 % spec_supercell.eps: 0x0 pixel, 300dpi, 0.00x0.00 cm, bb=14 220 581 621
\caption{Decay rates of $\log_{10}\left(\frac{|\lambda_L -  \lambda_{L_{\rm ref}}|}{\lambda_{L_{\rm ref}}}\right)$ (Eigenvalue), $\log_{10}\left(\frac{\|\phi_L -  \phi_{L_{\rm ref}}\|^2_{L^2(\R)}}{\|\phi_{L_{\rm ref}}\|_{L^2(\R)}^2}\right)$ (Error L2) and 
$\log_{10}\left(\frac{\|\phi_L -  \phi_{L_{\rm ref}}\|^2_{H^1(\R)}}{\|\phi_{L_{\rm ref}}\|_{H^1(\R)}^2}\right)$ (Error H1) for different values of $L$.}
\end{figure}

\medskip

The last series of numerical tests aims at testing the effect of numerical integration. For all $L\in\N^*$, we denote by 
$\lambda_{L,N_L,M_L}$ the eigenvalue of $H_{L,N_L,M_L}$ which is closest to $\lambda$, by $u_{L,N_L,M_L}$ an associated normalized eigenvector, and by $\phi_{L,N_L,M_L} = \chi_L u_{L,N_L,M_L}$ (we choose the sign of $u_{L,N_L,M_L}$ in such a way that $\|\phi_{L,N_L,M_L} -\phi_L\|_{L^2(\R^d)} \simeq 0$). In the plots below are drawn the errors $|\lambda_{L,N_L,M_L} - \lambda_L|$, $\|\phi_{L,N_L,M_L} -\phi_L\|_{L^2(\R^d)}$ and 
$\|\phi_L - \phi_{L,N_L,M_L}\|_{H^1(\R^d)}$ for the following values: 
\begin{itemize}
 \item $L = 6,8,10,12,14,16,18$,
\item $N_L = N L$ where $N = 2,4,6,8,10,12,14$,
\item $M_L = M L$ where $M = 56, 112, 224, 448$,
\end{itemize}
as well as the results obtained with exact integration ($M=\infty$).

\begin{figure}\label{fig:fig3}
\centering
 \includegraphics[width=14cm]{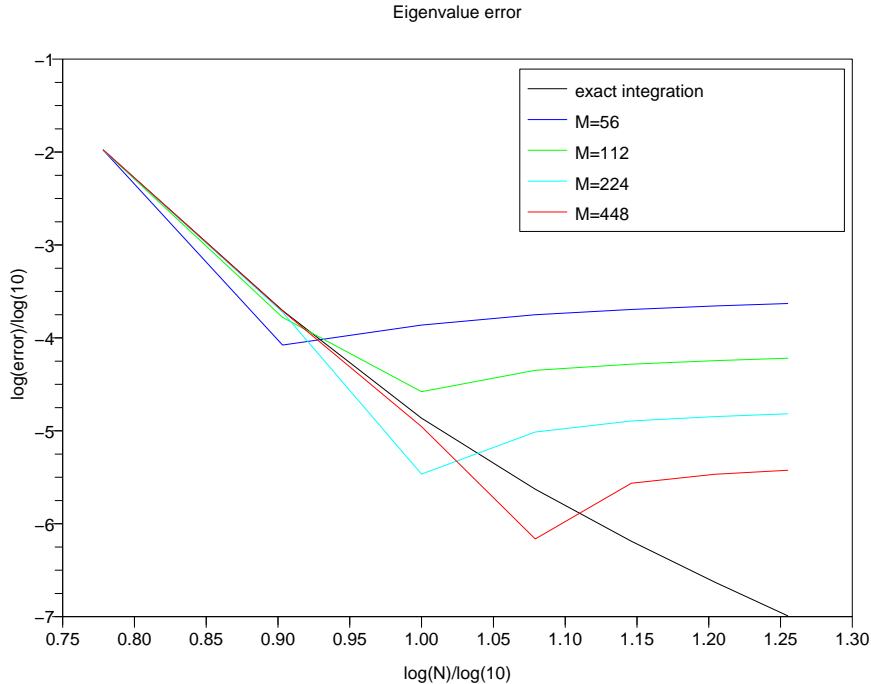}
 % spec_supercell.eps: 0x0 pixel, 300dpi, 0.00x0.00 cm, bb=14 220 581 621
\caption{Error on the eigenvalue $\log_{10}\left(|\lambda_{L,N_L,M_L} - \lambda_L|\right)$ as a function of $\log_{10}(N)$.}
\end{figure}

\begin{figure}\label{fig:fig4}
\centering
 \includegraphics[width=14cm]{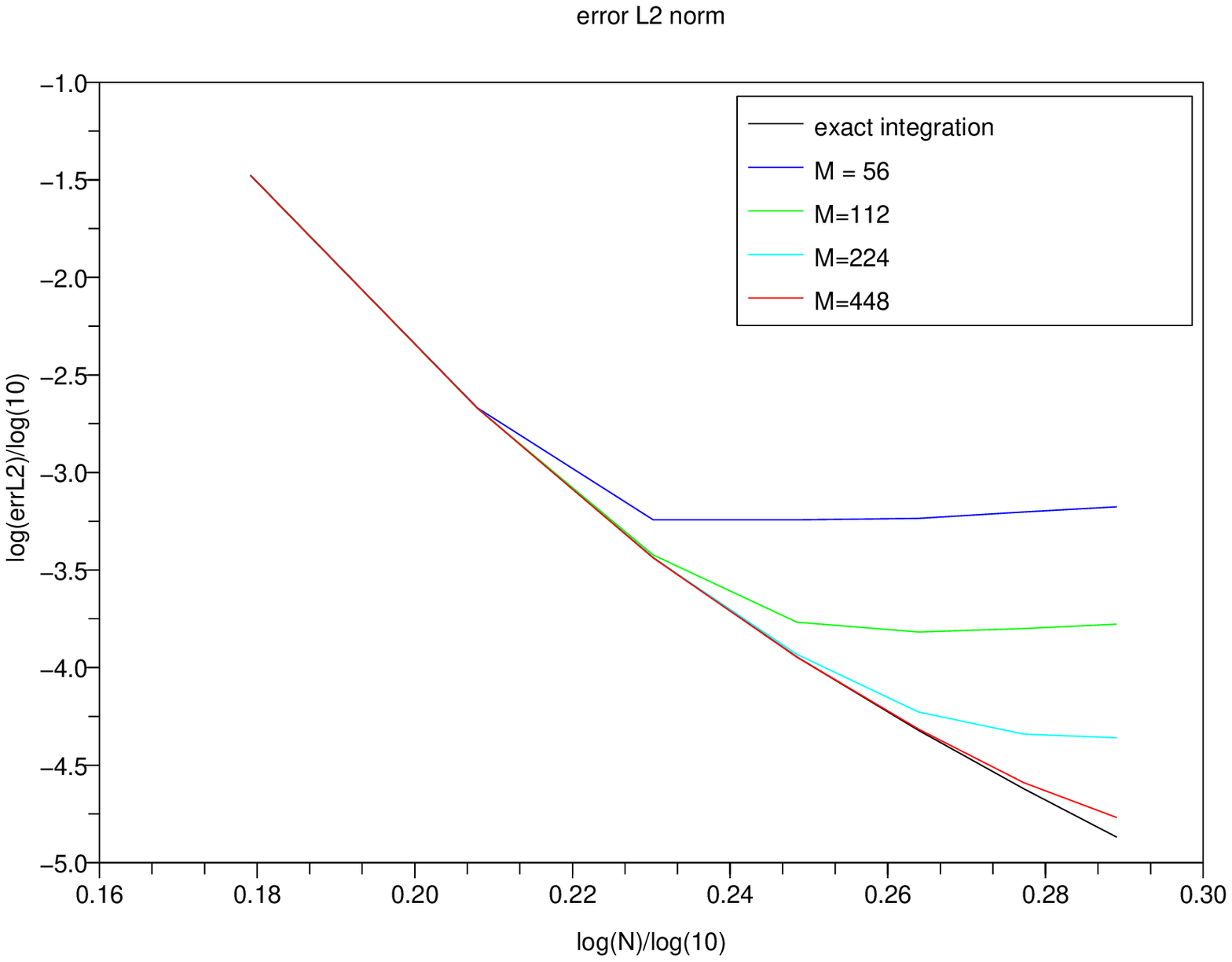}
 % spec_supercell.eps: 0x0 pixel, 300dpi, 0.00x0.00 cm, bb=14 220 581 621
\caption{Error on the eigenvector $\log_{10}\left(\|\phi_{L,N_L,M_L} - \phi_L\|_{L^2(\R^d)}\right)$ as a function of $\log_{10}(N)$.}
\end{figure}

\begin{figure}\label{fig:fig5}
\centering
 \includegraphics[width=14cm]{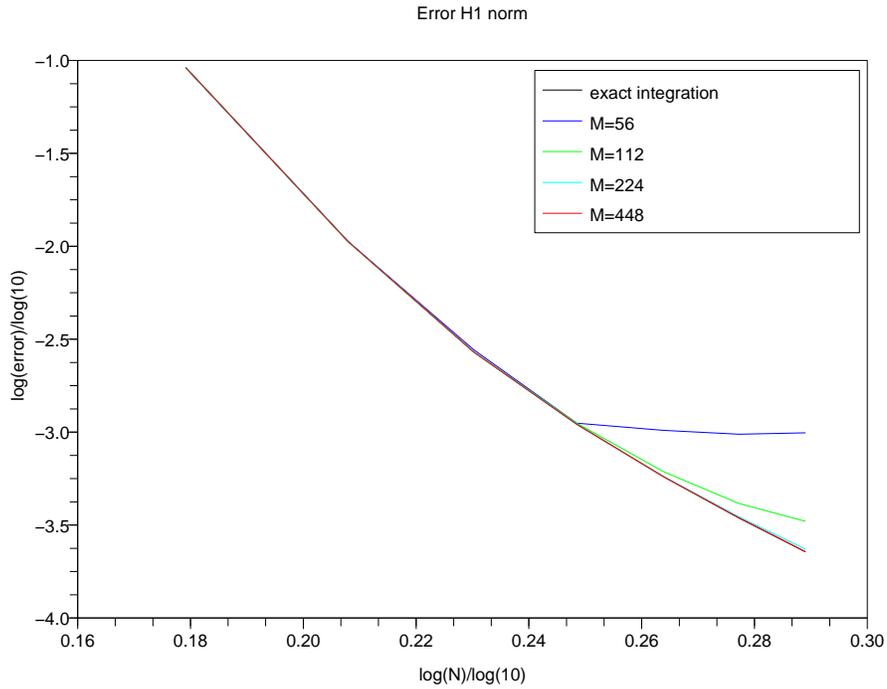}
 % spec_supercell.eps: 0x0 pixel, 300dpi, 0.00x0.00 cm, bb=14 220 581 621
\caption{Error on the eigenvector $\log_{10}\left(\|\phi_{L,N_L,M_L} - \phi_L\|_{H^1(\R^d)}\right)$ as a function of $\log_{10}(N)$. }
\end{figure}

\section{Appendix: Banach-Ne\v{c}as-Babu\v{s}ka's Theorem and Strang's lemma}\label{sec:appendix}

In this appendix, we recall the Banach-Ne\v{c}as-Babu\v{s}ka theorem and the Strang lemma (see e.g. \cite{Fortin,Ern}). 

\begin{theorem}\normalfont \bfseries (Banach-Ne\v{c}as-Babu\v{s}ka) \normalfont \itshape
Let $W$ be a Banach space and $V$ a reflexive Banach space. Let $a\in\cL(W\times V; \R)$ and $f\in V'$. Then the problem 
\begin{equation}
 \left\{
\begin{array}{l}
 \mbox{find }u\in W\mbox{ such that}\\
\forall v\in V, \quad a(u,v) = f(v),\\
\end{array}
\right .
\end{equation}
is well-posed if and only if
\begin{itemize}
 \item $\exists \alpha >0, \mbox{s.t.} \; \displaystyle \mathop{\inf}_{w\in W} \mathop{\sup}_{v\in V} \frac{|a(w,v)|}{\|w\|_W\|v\|_V}\geq \alpha$;
\item $\forall v\in V, \; \left( \forall w\in W, \quad a(w,v) = 0\right) \Rightarrow (v=0)$.
\end{itemize}
Moreover, the following a priori estimate holds:
\begin{equation}\label{eq:Banach}
 \forall f\in V', \quad \|u\|_W \leq \frac{1}{\alpha}\|f\|_{V'}.
\end{equation}
\end{theorem}

\begin{lemma}\normalfont \bfseries (Strang) \normalfont \itshape
Let us consider the following approximate problem
\begin{equation}
 \left\{
\begin{array}{l}
 \mbox{find }u_n\in W_n\mbox{ such that}\\
\forall v_n\in V_n, \quad a_n(u_n, v_n) = f_n(v_n),\\
\end{array}
\right .
\end{equation}
and let us assume that
\begin{itemize}
 \item $W_n \subset W$ and $V_n \subset V$;
\item $\exists \alpha_n>0,\mbox{s.t.}\;  \displaystyle \mathop{\inf}_{w_n\in W_n} \mathop{\sup}_{v_n\in V_n} \frac{|a_n(w_n, v_n)|}{\|w_n\|_W \|v_n\|_V}\geq \alpha_n$,
and $\mbox{\rm dim}(W_n) = \mbox{\rm dim}(V_n)$;
\item the bilinear form $a_n$ is bounded on $W_n\times V_n$.
\end{itemize}
Then, the following error estimate holds:
\begin{eqnarray*}
 \|u-u_n\|_W &\leq& \frac{1}{\alpha_n} \|f-f_n\|_{\cL(V_n)} \\
&& + \mathop{\inf}_{w_n\in W_n} \left[ \left( 1 + \frac{\|a\|_{\cL(W,V_n)}}{\alpha_n}\right)\|u-w_n\|_W + \frac{1}{\alpha_n}\mathop{\sup}_{v_n\in V_n} \frac{|a(w_n, v_n) - a_n(w_n, v_n)|}{\|v_n\|_V}\right].\\
\end{eqnarray*}
\end{lemma}

\bibliography{biblio}

\end{document}